\documentclass[twoside,a4paper,reqno]{amsart}

\usepackage{style}

\usepackage{graphicx}

\usepackage{pdfsync}

\usepackage{xcolor}
\definecolor{webgreen}{rgb}{0,.5,0}
\definecolor{webbrown}{rgb}{.6,0,0}
\definecolor{RoyalBlue}{cmyk}{1, 0.50, 0, 0}

\usepackage{mathrsfs}
\usepackage{courier} 
\normalfont

\usepackage{epsfig, graphicx, subfigure}
\usepackage{verbatim, setspace}
\usepackage[colorlinks=true, breaklinks=true, urlcolor=webbrown, linkcolor=RoyalBlue, citecolor=webgreen,backref=page]{hyperref}

\numberwithin{equation}{section}

\usepackage{courier} 
\normalfont

\let\Re=\undefined
\DeclareMathOperator{\Re}{Re}
\let\Im=\undefined
\DeclareMathOperator{\Im}{Im}

\def\ge{\geqslant}
\def\le{\leqslant}
\def\ol{\overline}

\def\norm[#1]{\left\| #1 \right\|}

\def\norm[#1]{\left\| #1 \right\|}
\def\snorm[#1]{\left\| #1 \right\|_{W}}

\theoremstyle{remark}

\begin{document}

\title[IST]{
Stability of Schur's iterates and fast solution \\ of the discrete integrable NLS
}

\author{R.\,V.\,Bessonov, \quad  P.\,V.\,Gubkin}

\address{
\begin{flushleft}
R.V.Bessonov: bessonov@pdmi.ras.ru\\\vspace{0.1cm}
St.\,Petersburg State University\\  
Universitetskaya nab. 7-9, 199034 St.\,Petersburg, RUSSIA\\
\vspace{0.1cm}
St.\,Petersburg Department of Steklov Mathematical Institute\\ Russian Academy of Sciences\\
Fontanka 27, 191023 St.Petersburg,  RUSSIA\\
\end{flushleft}
}
\address{
\begin{flushleft}
P.V.Gubkin: 
gubkinpavel@pdmi.ras.ru\\\vspace{0.1cm}
St.\,Petersburg State University\\  
Universitetskaya nab. 7-9, 199034 St.\,Petersburg, RUSSIA\\
\vspace{0.1cm}
St.\,Petersburg Department of Steklov Mathematical Institute\\ Russian Academy of Sciences\\
Fontanka 27, 191023 St.Petersburg,  RUSSIA\\
\end{flushleft}
}

\thanks{
The work of RB in Sections 1,\,3,\,6 is supported by the Russian Science Foundation grant RScF 19-71-30002. The work of PG in Sections 2,\,4,\,5 is supported by Ministry of Science and Higher Education of the Russian Federation, agreement 075–15–2022–287 
and in part by the M\"obius Contest Foundation for Young Scientists. RB is a Young Russian Mathematics award winner and would like to thank its sponsors and jury.}

\begin{abstract}
We prove a sharp stability estimate for Schur iterates of contractive analytic functions in the open unit disk. 
We then apply this result in the setting of the inverse scattering approach and obtain a fast algorithm for solving the discrete integrable nonlinear Schr\"odinger equation (Ablowitz-Ladik equation) on the integer lattice, $\Z$. We also give a self-contained introduction to the theory of the nonlinear Fourier transform from the perspective of Schur functions and orthogonal polynomials on the unit circle. 
\end{abstract} 

\vspace{1cm}

\subjclass[2010]{35Q55, 42C05}
\keywords{Ablowitz-Ladik equation, Inverse scattering, Schur's algorithm, Szeg\H{o} class}

\maketitle
	
\section{Introduction}
\subsection{Schur's algorithm} The Schur class $\Sch(\D)$ in the open unit disk $\D = \{z \in \C: |z|<1\}$ of the complex plane, $\C$, consists of analytic functions $F$ in $\D$ such that $$\sup_{z \in \D}|F(z)| \le 1.$$ 
For $F\in\Sch(\D)$, we write $F\in \Sch_{*}(\D)$ if $F$ is not a finite Blaschke product. Take $F\in \Sch_*(\D)$, set $F_0 = F$, and define the sequence $\{F_n\}_{n \ge 0}$ using Schur's algorithm:
\begin{equation}\label{eq20}
zF_{n+1} = \frac{F_n - F_n(0)}{1 - \ov{F_n(0)}F_n}, \qquad n \ge 0.
\end{equation}
By construction and Schwarz lemma, the resulting functions $F_0$, $F_1$, $F_2$, $\ldots$ will belong to the class $\Sch_{*}(\D)$ as well. In the case where $F \in \Sch(\D) \setminus \Sch_{*}(\D)$ is a Blaschke product  of order $N \ge 0$, the same construction gives a finite sequence of Blaschke products $F_0$, $F_1$, \ldots $F_N$ of orders $N$, $N-1$, \ldots $0$, correspondingly. In particular, $F_N$ is a constant of unit modulus and the Schur's algorithm stops. 
	
\medskip
	
Note that $|F_n(0)| < 1$ for each $F \in \Sch_*(\D)$, $n \ge 0$, by the maximum modulus principle.
Therefore, each function $F \in \Sch_*(\D)$ generates a sequence of numbers $\{F_n(0)\}_{n \ge 0} \subset \D$. They are called the recurrence coefficients of $F$. It can be shown that the mapping
$$
F \mapsto \{F_n(0)\}_{n \ge 0}
$$
is a homeomorphism from $\Sch_*(\D)$ with the topology of convergence on compact subsets of $\D$ onto the space of sequences $q: \Z_+\to \D$  with the topology of elementwise convergence, see Section 1.3.6 in \cite{Simonbook1}. Here, $\Z_+ = \Z \cap [0, +\infty)$. In particular, for every sequence $\{\alpha_n\}_{n \ge 0} \subset \D$ there exists a unique function $F \in \Sch_*(\D)$ such that $\alpha_n = F_n(0)$ for every $n \in \Z_+$. In this paper we study stability of Schur's algorithm. We prove a sharp estimate for $|F_n(0)-  G_n(0)|$ in terms of $F-G$ for functions $F, G \in \Sch_*(\D)$ from the  Szeg\H{o} class, whose definition we now recall.
	
\medskip
	
Let $m$ denote the Lebesgue measure on the unit circle $\T = \{z \in \C: |z|=1\}$ normalized by $m(\T) = 1$. The following theorem can be found, e.g., in Section 2.7.8 of \cite{Simonbook1}.
	
\begin{Thm}[Szeg\H{o} theorem]\label{szthm} Let $F \in \Sch_*(
\D)$, and let $\{F_n(0)\}_{n\ge 0}$ be its recurrence coefficients. Then 
$$
\int_{\T}\log(1-|F|^2)\,dm = \log\prod_{n=0}^{\infty}(1-|F_n(0)|^2),
$$
where both sides are finite or infinite simultaneously.
\end{Thm}
\noindent We will refer to functions $F \in \Sch_*(\D)$ such that 
\begin{gather}\label{eq: eta def}
	\eta(F) = \prod_{n=0}^{\infty}(1-|F_n(0)|^2) > 0,
\end{gather}
as Schur functions of Szeg\H{o} class. Given any $r \in (0,1)$, and an analytic function $F$ in $\D$, we set
$$
\|F\|_{L^2(r\T)} = \left(\int_{\T}|F(r\xi)|^2\,dm(\xi)\right)^{1/2}.
$$
\begin{Thm}\label{t0}
Let $F$, $G$ be Schur functions of Szeg\H{o} class, and let $\eta > 0$ satisfy $\min(\eta(F), \eta(G)) \ge \eta$. Then for every $r\in (0,1)$ and $n\in \Z_+$, the estimate  
\begin{equation}\label{t0est}
\|F_n - G_n\|_{L^2(r\T)}\le C(\eta, r)r^{-n}\|F - G\|_{L^2(r\T)} ,
\end{equation}
holds with the constant $C(\eta, r) = \exp\bigl( \log\eta^{-1}\cdot \bigl(2 + \frac{1}{1-\sqrt{1-\eta}}\bigr)(\frac{4}{(1 - r)^2}  + 1)\bigr)$ depending only on $\eta$, $r$.
\end{Thm}
The order of the exponential factor $r^{-n}$ in Theorem \ref{t0} is sharp. Indeed, one can take $\delta \in (0,1)$ and set $F =\delta z^n$, $G = 0$. Then $F_{n}(z) = \delta$, $G_n(z) = 0$ for all $z \in \D$. So, we have $\|F_n - G_n\|_{L^2(r\T)} = \|\delta\|_{L^2(r\T)} = \delta$ and $\|F - G\|_{L^2(r\T)} = \|\delta z^n\|_{L^2(r\T)} = \delta r^n$ in this case. Since $\eta(F) = 1-\delta^2$, $\eta(G) = 1$ do not depend on $n$, a consideration of large $n$'s shows that the order of growth $r^{-n}$ in \eqref{t0est} cannot be improved within the Szeg\H{o} class. 
	
	\medskip
	
Theorem \ref{t0} can be used to estimate $|F_n(0) - G_n(0)|$ if we know that Schur functions $F$, $G$ are sufficiently close to each other in the disk $|z|\le r$. Indeed, by Bessel inequality, we have
\begin{equation}\label{eq2}
|F_n(0) - G_n(0)| \le \|F_n - G_n\|_{L^2(r\T)},
\end{equation}
because the system $\{z^{k}\}_{k \ge 0}$ is orthogonal in $L^2(r\T)$. We want to emphasize that the constant $C(\eta, r)$ in Theorem \ref{t0} is uniform for functions $F\in \Sch_{*}(\D)$ with the Szeg\H{o} constant $\eta(F)$ separated from zero. This is the most important feature of \eqref{t0est} when it compared with another stability result from the inverse spectral theory -- Sylvester-Winebrenner theorem \cite{Sylvester}. In the language of Schur functions, this theorem says that Schur's algorithm defines a homeomorphism in appropriate metric spaces:
\begin{Thm}[Sylvester--Winebrenner theorem]\label{SWt} The mapping 
$F \mapsto \{F_n(0)\}_{n \ge 0}$ that takes a Schur function into the sequence of its recurrence coefficients is a homeomorphism from the metric space $X_+ = \{F \in \Sch_{*}(\D): \eta(F) >0 \}$ 
with the metric $\rho_s(F,G)^2 = -\int_{\T}\log\bigl(1 - \bigl|\frac{F-G}{1-\bar F G}\bigr|^2\bigr)\,dm$ onto the metric space $\ell^2(\Z_+, \D)$ of square summable sequences $q: \Z_+ \to \D$ with the metric 
$\|q - \tilde q\|_{\ell^2}^{2} = \sum_{n \in \Z_+}|q(n) - \tilde q(n)|^2$. 
\end{Thm}
We prove this version of  Sylvester--Winebrenner theorem in Section \ref{s6}. It is very natural to expect that the modulus of continuity of the homeomorphism in Theorem \ref{SWt} is controlled by $\eta$ on the subset of $F \in X_+$ with $\eta(F)> \eta > 0$. This is, however, not the case! See Proposition \ref{pr12} below. On the other hand, the uniform character of  estimate \eqref{t0est} will be crucial for the application of \eqref{t0est} to the discrete integrable nonlinear Schr\"odinger equation (Ablowitz-Ladik equation). Let us discuss it next.


\medskip

\subsection{AL: statement of the problem} \noindent Consider the defocusing Ablowitz-Ladik equation (AL) on the integer lattice, $\Z$,
\begin{equation}\label{al}
\frac{\partial}{\partial t}q(t,n) = i\bigl(1 - |q(t,n)|^2\bigr)\bigl(q(t,n-1) + q(t,n+1)\bigr), \quad q(0, n) = q_0(n), \quad n \in \Z. 
\end{equation}
The variable $t \in \R$ is considered as time,  $n \in \Z$ is the discrete space variable. Ablowitz-Ladik equation is the integrable model introduced in \cite{AL75}, \cite{AL76} as a spatial discretization of the cubic non-linear Schr\"odinger equation (NLS), see \cite{APT04} for a general context and modern exposition. If we change variables to $u = e^{-2it}q$, then \eqref{al} becomes
\begin{align*}
i\frac{\partial}{\partial t}u(t,n) 
&= -\bigl(1 - |u(t,n)|^2\bigr)\bigl(u(t,n-1) + u(t,n+1)\bigr) + 2u(t,n)  \\
&=-u(t,n-1) + 2u(t,n) - u(t,n+1) + |u(t,n)|^2(u(t,n-1) + u(t,n+1)),
\end{align*}
which is indeed a discretization of the continuous defocusing NLS equation,
\begin{equation}\label{nls}
i \frac{\partial}{\partial t}u(t,x) = -\frac{\partial^2}{\partial x^2}u(t,x) + 2|u(t,x)|^2u(t,x), \qquad x \in \R.
\end{equation}
We are going to present a new solution method for \eqref{al}  based on Schur's algorithm. The rate of its convergence will be estimated using Theorem \ref{t0}. We deal with the following problem:
	
	\medskip

\begin{Prob}\label{Problem}
Given $\eps \in (0,1)$, $t \in \R$, $n_0 \in \Z$, and a sequence $q_0$ on $\Z$ such that $|q_0(n)| < 1$, $n \in \Z$, 
$$\prod_{n \in \Z}(1-|q_0(n)|^2) \ge \eta > 0,$$ evaluate the solution $q$ of \eqref{al} at $(t,n_0)$ with the absolute error at most $\eps$.
\end{Prob}

\medskip

The quantity $\prod_{n \in \Z}(1-|q(t,n)|^2) = \prod_{n \in \Z}(1-|q_0(n)|^2)$ is conserved under the flow of AL equation. So, it is a natural characteristic for results on stability/accuracy of solutions of AL equation.

\medskip

We introduce the algorithm which solves Problem \ref{Problem} in $O(\mathbf{n}\log^{2}\mathbf{n})$ operations, where $\mathbf{n} = t + \log\eps^{-1}$. Thus, to have accuracy $e^{-n}$ at the moment of time $t = 1$, one need to take at most $c_{\eta} n \log^{2}n$ arithmetic operations for some constant $c_{\eta}>0$ depending only on $\eta$. The basic Runge-Kutta scheme RK4 requires $n$ time steps ($\sim n \cdot k$ operations) for computing $u(1, j)$, $-k \le j \le k$, to guarantee accuracy $O(1/n^4)$ if we additionally assume that the impact of $u(t,j)$, $|j| \ge k$, is negligible for $0 \le t \le 1$. 

\subsection{AL: localization} Our solution method is a modification of the classical inverse scattering approach. From a bird-view look, the standard procedure (see Chapter 2 in \cite{TTT}) of solving \eqref{al} by means of the inverse scattering theory (IST) looks as follows: given an initial datum $q_0: \Z \to \D$, define the so-called reflection coefficient $\rc_{q_0}$ by   
\begin{equation}\label{eq25bisbis}
	\rc_{q_0}(z) = \frac{b(z)}{a(z)}, \qquad 
	\begin{pmatrix}
		a & b\\
		\ov{b} & \bar a
	\end{pmatrix} 
	=
	\prod_{k \in \Z}
	\frac{1}{\sqrt{1-|q_0(k)|^2}}
	\cdot
	\begin{pmatrix}
		1 & \ov{q_0(k)} z^{-k}\\
		q_0(k)  z^k & 1
	\end{pmatrix}, \qquad z \in \T,
\end{equation}
and find $q(t,\cdot): \Z \to \D$ such that $\rc_{q(t, \cdot)} = e^{it(z + 1/z)}\rc_{q_0}$ on $\T$. It turns out that $q(t, \cdot)$ will solve \eqref{al} for initial datum $q_0$ provided $q_0$ decays fast enough (say, $\sum_{k \in \Z}|q_0(k)| < \infty$). A fundamental problem appearing when one tries to solve \eqref{al} by IST with merely $\ell^2(\Z, \D)$ initial datum $q_0$ (i.e., for general $q_0: \Z \to \D$ such that $\prod_{n\in \Z}(1 - |q_0(n)|^2) > 0$) is that we can have $\rc_{q_{0}} = \rc_{\tilde q_{0}}$ for $q_{0} \neq \tilde q_{0}$. This phenomenon was first observed by Volberg and Yuditskii in \cite{VYu02} on the level of Jacobi matrices, and then by Tao and Thiele \cite{TT} in the setting of the nonlinear Fourier transform, NLFT. It shows that when we pass to reflection coefficients $\rc_{q_{0}}$, $\rc_{\tilde q_{0}}$, some information gets lost and there are no chances to solve \eqref{al} for $\ell^2(\Z, \D)$ initial data by using IST approach directly. To overcame this difficulty (non-injectivity of NLFT), we first prove the following localization estimate.
\begin{Thm}\label{t2}
Let $q_0:\Z \to \D$ be such that $\prod_{n\in\Z}(1 - |q_0(n)|^2)\ge \eta$ for some $\eta > 0$ and let $q$ be the solution of \eqref{al} for the initial datum $q_0$. Take $N \in \Z_+$, consider the sequence $q_{0, N}$ defined by
\begin{gather*}
q_{0, N}(n) = 
\begin{cases}
q_{0}(n), &|n|\le N,\\
0, & |n|> N,
\end{cases}
\end{gather*}
and let $q_{N}$ be the corresponding solution of \eqref{al}. Then, for $N\ge |j|$, $t>0$ and all $r \in (0,1)$, we have
\begin{gather}\label{eq: theorem 1.3 ineq}
|q(t, j) - q_N(t,j)|\le \frac{4e^{t/r}C(\eta, r)}{1 - r} r^{N-|j|},
\end{gather}
where $C(\eta, r)$ is the function from Theorem \ref{t0}.
\end{Thm}
Having in mind a possible future development of a parallel theory for continuous NLS equation, we use only ``spectral'' methods in the proof of Theorem \ref{t2}. The reader interested in short and elementary proof of Theorem \ref{t2} by means of a direct approach, could find it in Section \ref{appendix}. 

\subsection{AL: compactly supported initial data}\label{section: Solution for compactly supported initial data} Having Theorem \ref{t2}, it remains to solve \eqref{al} for compactly supported initial data $q_{0}:\Z \to \D$. This can be done by a variety of methods, both theoretically and numerically. In particular, the standard IST approach works in this case, but accuracy estimates for numerical schemes based on IST and $\ell^2$-bounds are missed in the literature. Taking into account the non-injectivity of NLFT, we see that the problem, in fact, is fairly nontrivial: some distant compactly supported data $q_0$, $\tilde q_0$ correspond to almost identical reflection coefficients $\rc_{q_0}$, $\rc_{\tilde q_0}$. Indeed, it is enough to take different $q_0, \tilde q_0 \in \ell^2(\Z, \D)$ with the same reflection coefficient and consider  restrictions of $q_0$, $\tilde q_0$ to a large discrete interval $[-N, N]$. Then the corresponding reflection coefficients will almost coincide by continuity of NLFT. This phenomenon, when ignored,  leads to instabilities. Below we describe a procedure that can be used to get the solution with prescribed accuracy. 

\medskip

Consider $q_0: \Z \to \D$ supported on $\Z\cap [0, \ell]$ for some $\ell \in \Z_+$. Note that $q$ solves \eqref{al} if and only if $q(t,\cdot + j)$ solves \eqref{al} for the initial datum $q_0(\cdot + j)$. Therefore, we do not loss generality when assuming $\supp q_0 \subset [0, \ell]$. Moreover, it is easy to see that $q(t, n)$ solves \eqref{al} if and only if $- q(-t, n)$ solves \eqref{al} with the initial data $-q_0$. So, we can also assume that $t>0$. 

\medskip

\noindent Consider the Fourier expansion of the inverse scattering multiplier $e^{it(z + 1/z)}$:
\begin{gather}\label{bessel generating function}
	e^{it(z + 1/z)} =\sum_{k \in \Z} i^k J_k(2t) z^k, \qquad z \in \T.
\end{gather}
Here, $J_k$ are the standard Bessel functions \cite{Bowman} of order $k$, i.e.,
\begin{gather*}
	J_k(2t) = i^{-k}\int_{\T}e^{it(z+1/z)}\bar z^{k}\,dm = \sum_{m = 0}^{\infty}\frac{(-1)^m t^{2m + k}}{(m + k )!}.
\end{gather*}
Let $P_{n, t} = \sum_{|k| \le n} i^k J_k(2t) z^k$ be the Laurent trigonometric polynomial of $e^{it(z + 1/z)}$ of order $n$. Define the function $G_{n, t}$ by
\begin{gather}\label{definition of G_n,t}
	G_{n, t} = (1 - \delta_{n,t})z^nP_{n,t},\quad \delta_{n,t} =  \frac{t^ne^t}{n!}.
\end{gather}
We will be interested in the situation when $n > ct$ with some $c > e$. In this case this ``$\delta_{n,t}$-correction'' is very small but important: it places $G_{n,t}$ into Schur class.
Given a sequence $q_0: \Z \to \D$ supported on $[0, \ell]$, define the  coefficients $a$, $b$ of $q_0$ by \eqref{eq25bisbis}. Note that the product in \eqref{eq25bisbis} contains at most $\ell+1$ nontrivial terms. One can check that $a$, $\bar b$ in \eqref{eq25bisbis} coincide on $\T$ with analytic polynomials in $z$ of degree at most $\ell$, and, moreover, $|\bar b(z)| < |a(z)|$ if $|z| \le 1$. Set $\fc_{q_0} = \bar b/a$. 
The function $F_{n,0}=G_{n,t}\fc_{q_0}$ is  rational and belongs to the Schur class $\Sch_*(\D)$ (see Proposition \ref{lem: modulus gnt} below). Fix $j \in \Z$ and use Schur's algorithm \eqref{eq20}
to find rational functions $F_{n,0}$, $F_{n,1}$, $F_{n,2}, \ldots, F_{n,n+j}, \ldots$ (Schur iterates of $F_{n, 0}$). Set
\begin{equation}\label{approx}
\tilde q_n(t, j) 
= 
\begin{cases}
F_{n,n+j}(0), & j \ge -n,\\
0 & j < -n.
\end{cases}
\end{equation}
The following theorem shows that $\tilde q_n$ approximates the solution $q$ of \eqref{al} with very high accuracy.
	
\begin{Thm}\label{t3}
Let $t > 0$, and let $q_0: \Z \to \D$ be a sequence compactly supported on $\Z_+$. Assume that $\prod_{n\in \Z_+}(1-|q_0(n)|^2) \ge \eta$ for some $\eta > 0$. Then, the function $\tilde q_n$  in \eqref{approx} satisfies
\begin{gather}\label{eq41}
|q(t, j) - \tilde q_n(t, j)| 	\le 2^{j}C(\eta, 1/2)\frac{12e^{5t}}{\sqrt{2\pi n}}\left(\frac{2et}{n}\right)^{n}, 
\end{gather}
for all $n \in \Z_+$, $j \in \Z$, $t > 0$ such that $n+j \ge 0$, $n > t$, and $\delta_{n,t} <1$, see \eqref{definition of G_n,t}. Here $q$ is the solution of \eqref{al} and $C(\eta, r)$ is the function from Theorem \ref{t0}.
\end{Thm}
Note that the right hand side in \eqref{eq41} is very small when $n$  is much larger than $2et$ and $j$ is fixed. The estimate does not depend on the size of the support of $q_0$. In fact, Theorem \ref{t3} remains true if we assume only $\supp q_0 \subset [0, +\infty)$ and $\prod_{0}^{+\infty}(1-|q_0(n)|^2) > 0$. In this case, it is known that the product in \eqref{eq25bisbis} converges in Lebesgue measure on $\T$ (see Section~\ref{s6}) and defines coefficients $a$, $b$ almost everywhere on $\T$. Moreover, $\fc_{q_0} = \bar b/a$ will coincide with non-tangential values of a function of Schur class $\Sch_*(\D)$. Then $\tilde q_n(t, j)$ are well-defined by \eqref{approx}, and \eqref{eq41} will hold for them. 
	
\subsection{AL: algorithm for Problem \ref{Problem}}\label{s1-4} Let us summarize the algorithm that solves Problem \ref{Problem} based on Theorems \ref{t2} and \ref{t3}. At first, one need to choose a window $\Delta = [n_0-N, n_0+N]$ where $N$ is such that $|q(t, n_0) - q_{N}(t,n_0)| \le \eps/2$ for the exact solution $q_{N}$ with the truncated initial datum $q_{0, N} = \chi_{\Delta} q_0$. Then, one need to shift $q_{0, N}$ by $n_0-N$ to make it supported on $\Z_+\cap [0, 2N]$ and use the algorithm described in Section \ref{section: Solution for compactly supported initial data} to find the approximate solution $\tilde q_n$ with accuracy $\eps/2$ at $j = N$ for the shifted sequence. Taking $N = 5 + [4et + \log_2 \frac{C(\eta, 1/2)}{\eps}]$, $n = 2N$, we will get $|\tilde q_n(t, N) - q(t, n_0)| \le \eps$, see Section~\ref{s7}. In Section~\ref{s7} we check that the whole procedure requires $O(\mathbf{n}\log^2\mathbf{n})$ operations for $\mathbf{n} = t + \log\eps^{-1}$. In fact, the sequence $\tilde q_n$ approximates $q$ with accuracy $O(\eps)$ on the interval $[n_0-N/2, n_0]$, not only at the point $n_0$. Considering reflection of $q_0$ and applying the algorithm twice, one can construct approximation to $q$ on $[n_0-N/2,n_0+N/2]$ in $O(\mathbf{n}\log^2\mathbf{n})$ operations.
	
	\medskip
	
\subsection{AL: historical remarks and motivation} 
As a classical integrable model, Ablowitz-Ladik equation has a well-developed theory in the periodic case \cite[Chapter 11]{Simonbook2}, \cite{N05}, \cite{MM23}, finite case \cite{GN09}, \cite{KN07},  
in the half-infinite case \cite{Gol06}, \cite{Sim07},  
and on the whole lattice $\Z$, 
see \cite{TTT}, \cite{GHMT08}, \cite{GHMT09}, \cite{KOVW23}. Paper \cite{GHMT09} contains a historical overview and an extensive 
bibliography, including works following original approach of Ablowitz and Ladik, who obtained a Lax pair for \eqref{al} by discretizing the Zakharov-Shabat Lax pair for the continuous NLS equation. Somewhat opposite, references mentioned in this paragraph (and results used in this paper) are mostly related to recent works that appeared after Nenciu and Simon \cite[Chapter 11]{Simonbook2}, \cite{N05} discovered a new Lax pair for this equation, making a connection to CMV matrices and orthogonal polynomials on the unit circle.
The IST method as a tool for existence theorems for Ablowitz-Ladik equation attracted a limited attention in the literature because the solvability of \eqref{al} for all initial data $q_0: \Z \to \D$ can be easily obtained by means of a fixed point theorem (see Appendix). 
However, Ablowitz-Ladik equation is a perfect model for developing an accurate fast IST-based numerical scheme that can be later generalized for the continuous NLS equation.

\medskip

\subsection{The nonlinear Fourier transform} The last part of the paper can be regarded as the introduction to the theory of the nonlinear Fourier transform, NLFT. The main results in this area are due to Thiele and Tao, see the paper \cite{TT} or its extended version by Thiele, Tao, and Tsai \cite{TTT}, where Ablowitz-Ladik equation appears in the setting of NLFT. Papers \cite{TT}, \cite{TTT} influenced much on the present work. We decided to give a short introduction to the theory of NLFT in the language of Schur functions and orthogonal polynomials to make the paper more self-contained. We hope that our arguments will be of independent interest for the orthogonal polynomials community.   

\medskip 

For $1 \le p < \infty$, let $\ell^p(\Z, \D)$ be the set of sequences $q: \mathbb \Z \to \D$ such that $\sum_{n \in \Z}|q(n)|^p < \infty$. We endow it with the usual distance $\|q_1 - q_2\|_{\ell^p} = \bigl(\sum_{n \in \Z}|q_1(n) - q_2(n)|^p\bigr)^{1/p}$. Note that $\ell^p(\Z, \D)$ is not a linear space. Using formula \eqref{eq25bisbis}, define the nonlinear Fourier transform (or the scattering map) by
$$
\F_{sc}: q \mapsto \rc_q,
$$
on the set $\ell^1(\Z, \D)$. Here we consider $\F_{sc}$ as the map from $\ell^1(\Z, \D)$ to $L^\infty(\T)$. Later on, the domain of $\F_{sc}$ will be extended, while the target space will be changed to a narrower one. Define the metric space
\begin{gather}\label{eq: X space def}
	X = \{h \in L^{\infty}(\T): \|h\|_{L^\infty(\T)} \le 1, \; \log(1-|h|^2) \in L^1(\T)\},
\end{gather}
with the Sylvester--Winebrenner metric $\rho_s$ (see \cite{Sylvester}) given by
\begin{gather}\label{eq: sylvester metric definition}
	\rho_s(h_1,h_2) = \sqrt{-\int_{\T}\log\left(1 - \Bigl|\frac{h_1-h_2}{1-\bar h_1 h_2}\Bigr|^2\right)\,dm}.
\end{gather}
For $\delta \in [0, 1)$, denote $B[\delta] = \{h \in L^\infty(\T):\; \|h\|_{L^\infty(\T)} \le \delta\}.$ We have $B[\delta] \subset X$ for every $\delta \in [0,1)$. So, let us consider $B[\delta]$ as the subspace of $X$ with induced metric topology. As we will see below, $\F_{sc}$ uniquely extends to the continuous map from $\ell^2(\Z, \D)$ to $X$. Set $\cG[\delta] = \F_{sc}^{-1}(B[\delta])$
where $\F_{sc}^{-1}(E)$ is the full preimage of a set $E$ under the mapping $\F_{sc}: \ell^2(\Z, \D) \to X$. 

\medskip

With this definitions at hand, we are ready to summarize the basic properties of $\F_{sc}$. 
\begin{Thm}\label{t5}
The nonlinear Fourier transform $\F_{sc}$ has the following properties:
\begin{itemize}
\item[$(1)$] the map $\F_{sc}: \ell^{1}(\Z, \D) \to L^\infty(\T)$ extends uniquely to the continuous map $\F_{sc}: \ell^2(\Z, \D) \to X$;
\item[$(2)$] the map $\F_{sc}: \ell^2(\Z, \D) \to X$ is closed; 
\item[$(3)$] we have $\F_{sc}(q(\cdot - n)) = z^{-n} \F_{sc}(q)$ for every $q \in \ell^2(\Z, \D)$, $n \in \Z$; 
\item[$(4)$] the map $\F_{sc}: \ell^2(\Z, \D) \to X$ is surjective; 
\item[$(5)$] the map $\F_{sc}: \ell^2(\Z, \D) \to X$ is not injective; 
\item[$(6)$] the map $\F_{sc}: \cG[\delta] \to B[\delta]$ is a homeomorphism for every $\delta \in (0, 1)$; 
\item[$(7)$] if $q = q(t, n)$ is the solution of \eqref{al} with the initial datum $q_0 \in \cG[\delta]$, then $q(t, \cdot) \in \cG[\delta]$ for each $t \in \R$, and 
$q(t, \cdot) = \F_{sc}^{-1}(e^{it(z+ 1/z)}\F_{sc}(q_0))$. 
\end{itemize}
\end{Thm}
Assertion $(2)$ in Theorem \ref{t5} is new. It implies, in particular, that $\F_{sc}$ is a homeomorphism on the set of potentials $q \in \ell^2(\Z, \D)$ that are completely determined by the reflection coefficient $\rc_{q}$. Assertion $(7)$ is not proved in our paper (we did not found a sufficiently short argument), see \cite{TTT} for the proof. Some ideas in the proof of Theorem \ref{t5} are due to S.\,Denisov, the authors would like to thank him for his contribution.  
	
\section{Schur's algorithm. Proof of Theorem \ref{t0}.}
In this section we prove Theorem \ref{t0}. For an analytic function $F$ in $\D$, it will be convenient to set
$$
M_{F}(r) = \max_{|z| \le r}|F(z)|.
$$
At first, we prove the following lemma.
\begin{Lem}\label{khr lemma}
Let $F \in \Sch_*(\D)$, and let $F_k$ be its Schur iterates defined by \eqref{eq20}. Then 
\begin{equation}\label{eq38}
\sum_{k = 0}^{\infty}M_{F_k}^2(r) \le \frac{4}{(1-r)^2}\cdot\sum_{k=0}^{\infty}|F_k(0)|^2 \le \frac{4}{(1-r)^2}\cdot\log(\eta(F)^{-1}), \qquad r \in [0, 1),
\end{equation}
where $\eta(F)$ is defined by \eqref{eq: eta def}.
\end{Lem}
\beginpf
Let us check the second inequality first. For $x\in (0,1)$, we have $(1 - x)^{-1}\ge e^x$, therefore
\begin{gather*}
	\eta(F)^{-1} = \prod_{k\ge 0}(1 - |F_k(0)|^2)^{-1}\ge \prod_{k\ge 0}e^{ |F_k(0)|^2} = e^{\sum_{k\ge 0}|F_k(0)|^2},
\end{gather*} 
which implies the required bound $\log(\eta(F)^{-1})\ge \sum_{k\ge 0}|F_k(0)|^2$. Now we focus on the first inequality in \eqref{eq38}. Set $\alpha_j = F_j(0)$, $j \ge 0$. We will use the estimate $(1.3.58)$ in \cite{Simonbook1} which reads
$$
|F(z)| \le 2\sum_{j = 0}^{\infty}|\alpha_{j}||z|^{j}, \qquad z\in \D.
$$
Applying it to $F_k$ in place of $F$ for $|z| = r$, we get $M_{F_k}(r) \le 2\sum_{j = 0}^{\infty}|\alpha_{k + j}|r^{j}$, hence
\begin{align*}
M_{F_k}^2(r)  \le 4\left(\sum_{j = 0}^{\infty}|\alpha_{k + j}|r^{j/2}\cdot  r^{j/2}\right)^2
\le 4\sum_{j = 0}^{\infty}|\alpha_{k + j}|^2r^{j}\cdot\sum_{j=0}^{\infty} r^{j} = \frac{4}{1 - r}\sum_{j = 0}^{\infty}|\alpha_{k + j}|^2r^{j}, 
\end{align*}
by Cauchy inequality.
Summing up over $k \in \Z_+$, we get
\begin{align*}
\sum_{k = 0}^{\infty}M_{F_k}^2(r) 
= \frac{4}{1 - r}\sum_{k = 0}^{\infty}\sum_{j = 0}^{\infty}|\alpha_{k + j}|^2r^{j}
&= \frac{4}{1 - r}\sum_{s = 0}^{\infty}|\alpha_s|^2\sum_{j = 0}^{s}r^{j} \le \frac{4}{(1 - r)^2}\sum_{s = 0}^{\infty}|\alpha_s|^2.
\end{align*}
This ends the proof. \qed 
\begin{Rema}Lemma \ref{khr lemma} holds with a better (for small $r$) estimate with $\frac{1}{(1 - r)^{4}}$ in place of $\frac{4}{(1 - r)^2}$. To prove this, one need to use expression for $F_k$ from Theorem 8.70 in \cite{Khr08}. A consideration of functions $F = \delta z^n$  for large $n$'s and small $\delta$'s shows that the constant in Lemma \ref{khr lemma} cannot be smaller than $\frac{1}{1-r^2}$.
\end{Rema}
	
	\medskip
	
\noindent{\bf Proof of Theorem \ref{t0}.}
Let $F, G \in \Sch_*(\D)$. We have
\begin{gather}\label{one step of schur algo}
	z(F_1 - G_1) 	= \frac{F - F(0)}{1 - \ol{F(0)}F} - \frac{G - G(0)}{1 - \ol{G(0)}G} = \frac{P}{Q}.
\end{gather}
Here, the numerator is 
\begin{align*}
	P &= (F - F(0))(1 - \ol{G(0)}G) - (G - G(0))(1 - \ol{F(0)}F) 
	\\
	&= \bigl[F - F(0) - G + G(0)\bigr] + FG(\ol{F(0)} - \ol{G(0)}) + \bigl[F(0)\ol{G(0)}G - \ol{F(0)}G(0)F\bigr].
\end{align*}
We have 
\begin{align*}
F(0)\ol{G(0)}G - \ol{F(0)}G(0)F 
=& F(0)\ol{G(0)}(G - F) + FF(0)(\ol{G(0)} - \ol{F(0)}) + \\
&+F\ol{F(0)} (F(0) - G(0)).
\end{align*}
It follows that 
\begin{align*}
	\|P\|_{L^2(r\T)} \le& \|F - F(0) - G + G(0)\|_{L^2(r\T)} + M_F(r)M_G(r) |F(0) - G(0)|+
	\\
	&+ |F(0)||G(0)|\|F - G\|_{L^2(r\T)} + 2M_F(r)|F(0)||F(0) - G(0)|.
\end{align*}
For an analytic function $H$ in $\D$, we have
\begin{gather*}
	|H(0)|\le \|H\|_{L^2(r\T)}, \quad \|H - H(0)\|_{L^2(r\T)}\le \|H\|_{L^2(r\T)},
\end{gather*}
by orthogonality of system $\{z^k\}_{k \ge 0}$. Applying this to $H = F- G$ and using $2xy \le x^2 + y^2$, we get
\begin{align*}
	\|P\|_{L^2(r\T)} \le& \|F - G\|_{L^2(r\T)}\left(1 + M_F(r)M_G(r) + |F(0)||G(0)| + 2M_F(r)|F(0)|\right),
	\\
	\le& \|F - G\|_{L^2(r\T)}\left(1 +\frac{3M_F(r)^2 + M_G(r)^2 + 3|F(0)|^2 + |G(0)|}{2} \right).
\end{align*}
Since $|P|$ remains the same when we swap $F$, $G$, we also have 
$$
\|P\|_{L^2(r\T)} \le \|F - G\|_{L^2(r\T)}\left(1 +\frac{M_F(r)^2 + 3M_G(r)^2 + |F(0)|^2 + 3|G(0)|}{2} \right).
$$
Taking a half-sum, we get  
\begin{equation}\label{eq12}
	\|P\|_{L^2(r\T)} \le \|F - G\|_{L^2(r\T)}\left(1 + M_F^2(r) + M_G^2(r) + |F(0)|^2 + |G(0)|^2\right).
\end{equation}
Further, for $z \in r\T$, we estimate the denominator $Q$ in \eqref{one step of schur algo} as follows: 
\begin{gather*}
	|Q(z)| = |(1 - \ol{G(0)}G)(1 - \ol{F(0)}F)|\ge (1 - |G(0)|M_G(r))(1 - |F(0)|M_F(r)),
\end{gather*}
where we use the fact that both brackets above are positive. Substitution of the bounds for $P$, $Q$ into \eqref{one step of schur algo} gives
\begin{gather*}
	r\|F_1 - G_1\|_{L^2(r\T)}\le \|F - G\|_{L^2(r\T)}\frac{1 + M_F^2(r) + M_G^2(r) + |F(0)|^2 + |G(0)|^2}{(1 - |G(0)|M_G(r))(1 - |F(0)|M_F(r))}.
\end{gather*}
The latter inequality applied to $F_k$ and $G_k$ in place of $F$, $G$ for $k = 0, \ldots, n - 1$ implies
\begin{gather}\label{eq: product of Ck}
	r^n \|F_n - G_n\|_{L^2(r\T)}\le \|F - G\|_{L^2(r\T)}\prod_{k = 0}^{n - 1} C_k,   
\end{gather}
for 
$$
C_k= \frac{1 + M_{F_k}^2(r) + M_{G_k}^2(r) + |F_k(0)|^2 + |G_k(0)|^2}{(1 - |G_k(0)|M_{G_k}(r))(1 - |F_k(0)|M_{F_k}(r))}.
$$
It remains to estimate $\prod_{k = 0}^{n-1} C_k$. For $\delta  \in (0, 1)$, denote by $c(\delta)$ the minimal positive number such that $\frac{1}{1 - x} \le 1 + c(\delta) x$ for all $x \in (0,1)$ satisfying $1-x^2>\delta$.
It is not difficult to check that
\begin{gather}\label{eq: c(delta) fromula}
	c(\delta) = \frac{1}{1-\sqrt{1-\delta}} \in \left[1,\; \frac{2}{\delta}\right].
\end{gather}
Observe that
\begin{gather}\label{eq: denom 1 - F_k^2 bound}
	1-|F_k(0)|^2 M_{F_k}(r)^2 \ge 1-|F_{k}(0)|^2 \ge \prod_{m=0}^{+\infty}(1-|F_{m}(0)|^2) = \eta(F) > \eta
\end{gather}
by our assumption. Then,
$$
\frac{1}{1-|F_k(0)|M_{F_k}(r)} \le 1+c(\eta)|F_k(0)|M_{F_k}(r) \le 1+\frac{c(\eta)}{2}(M_{F_k}^2(r) + |F_k(0)|^2).
$$
A similar estimate holds for functions $G_k$. It follows that 
\begin{align*}
	C_{k} \le& (1 + M_{F_k}^2(r) + M_{G_k}^2(r) + |F_k(0)|^2 + |G_k(0)|^2)
	\\
	\times&\left(1+\frac{c(\eta)}{2}(M_{F_k}^2(r) + |F_k(0)|^2)\right)\left(1 + \frac{c(\eta)}{2}(M_{G_k}^2(r) + |G_k(0)|^2)\right) 
	\\
	\le& \exp\left( \left(1 + \frac{c(\eta)}{2}\right)(M_{F_k}^2(r) + M_{G_k}^2(r) + |F_k(0)|^2 + |G_k(0)|^2)\right), 
\end{align*}
where we used the elementary inequality $1+x \le e^{x}$ three times.
Then, from Lemma \ref{khr lemma} we get
\begin{align*}
\prod_{k = 0}^{n-1}C_k \le \exp\left(\left(1 + \frac{c(\eta)}{2}\right)\left(\sum_{k = 0}^{n-1} M_{F_k}^2(r) + \sum_{k = 0}^{n-1} M_{G_k}^2(r) + \sum_{k = 0}^{n-1} |F_k(0)|^2 + \sum_{k = 0}^{n-1}|G_k(0)|^2\right)\right)
\\
\le \exp\left(\left(1 + \frac{c(\eta)}{2}\right)\left(\frac{8\log\eta^{-1}}{(1 - r)^2}  + 2\log\eta^{-1}\right)\right). 
\end{align*}
Substitution of the latter into \eqref{eq: product of Ck} and the bound \eqref{eq: c(delta) fromula} imply \eqref{t0est} with 
\begin{gather}\label{eq: t0 constant}
	C(\eta, r) = \exp\left( \log\eta^{-1}\cdot \left(2 + \frac{1}{1-\sqrt{1-\eta}}\right)\left(\frac{4}{(1 - r)^2}  + 1\right)\right).
\end{gather}
This ends the proof. \qed	
\begin{Rema}
The function $C(\eta, r)$ is very large if $\eta$ is not close to $1$ or if $r$ is close to $1$.
We have, e.g., $5 \cdot 10^{27} \le C(1/2, 1/2) \le  6\cdot 10^{27}$, $10^6 \le C(4/5, 1/2) \le  2\cdot 10^6$, and $9 \le C(24/25, 1/2) \le 10$.
 In \cite{KOVW23}, Killip, Ouyang, Visan, and
Wu proved that the continuous NLS equation with arbitrary $L^2(\R)$-initial data can be approximated by the solutions of equation \eqref{al}. It is interesting to note that $\eta \to 1$ in their construction during approximation process.    
\end{Rema}

\medskip
	
\section{Estimates for the multipliers. Proof of Theorem \ref{t3}}\label{s3}
Recall the definition \eqref{definition of G_n,t} of $G_{n,t}$ and $P_{n,t}$: 
\begin{gather}\label{definition of G_n,tbis}
P_{n, t} = \sum_{|k| \le n} i^k J_k(2t) z^k, \qquad 	G_{n, t} = (1 - \delta_{n,t})z^nP_{n,t},\qquad \delta_{n,t} =  \frac{t^ne^t}{n!}, \qquad t> 0.
\end{gather}
In this section we first prove a bound for $G_{n,t}$ and estimate the rate of convergence of $G_{n + 1,t} - zG_{n,t}$ to zero. Then we prove Theorem \ref{t3}. Throughout this section, we assume that $t>0$. 
\begin{Lem}\label{lem: modulus gnt}
Let $z\in\T$ and let $n > t > 0$ be such that $\delta_{n, t} < 1$ for $\delta_{n,t} = \frac{t^ne^t}{n!}$ from \eqref{definition of G_n,t}. Then we have $|G_{n,t}(z)| < 1$. In particular, for every $q_0 \in \ell^2(\Z, \D)$ with $\supp q_0 \subset \Z_+$, we have $G_{n,t}\fc_{q_0} \in \Sch_*(\D)$ and the construction described in Section \ref{section: Solution for compactly supported initial data} is correct.
\end{Lem}
\beginpf
We have
\begin{gather*}
\left|P_{n,t}(z) - e^{it(z + 1/z)}\right| = \left|\sum_{|k| > n} i^k J_k(2t) z^k\right| \le 2\sum_{k > n} r^{-k}|J_k(2t)|, \qquad |z| = r.
\end{gather*}
The standard estimate (see, e.g., page 91 in \cite{Bowman}) $|J_{\nu}(2t)|\le {|t|^{\nu}/\Gamma(\nu + 1)}$ implies
\begin{align}\label{p and exp diff}
\left|P_{n,t}(z) - e^{it(z + 1/z)}\right| &\le 2\sum_{k > n} \frac{r^{-k}t^{k}}{k!} \le \frac{2t^nr^{-n} e^{t/r}}{(n+1)!} \le \frac{t^nr^{-n} e^{t/r}}{n!},
\\
\label{pnt bound}
\left|P_{n,t}(z)\right| &\le\left|e^{it(z + 1/z)}\right| + \frac{t^nr^{-n} e^{t/r}}{n!}.
\end{align}
In particular, for $z\in \T$ this gives $|P_{n,t}(z)|\le 1 + \delta_{n,t}$, where $\delta_{n,t} = \frac{t^ne^t}{n!}$ is from \eqref{definition of G_n,t}. Therefore, we have
\begin{gather*}
|G_{n,t}(z)| = (1 - \delta_{n,t})|P_{n,t}(z)| \le 1 - \delta_{n,t}^2 < 1, \qquad z \in \T,
\end{gather*}
where the factor $(1 - \delta_{n,t})$ is positive by our assumption. For compactly supported $q_0$ with $\supp q_0 \subset [0, \ell]$, it is not difficult to check that $\fc_{q_0}$ is a Schur function by considering partial products in \eqref{eq25bisbis} and using induction. For the general case, see formula \eqref{eq46} below. Then, we have  $G_{n,t}\fc_{q_0} \in \Sch_*(\D)$ by construction. \qed
\begin{Lem}\label{Lemma: convergence of gnt}
Let $n$, $t$ be as in Lemma \ref{lem: modulus gnt}. Then we have $\max_{|z|=r}|G_{n, t}| \le  e^{t/r}(r^n + 3\delta_{n,t})$ for $r\in(0,1)$, and, moreover,
\begin{gather*}
\max_{|z|=r}\left|G_{n + 1,t}(z) - zG_{n,t}(z)\right| \le  S_n(t,r), \quad S_n(t,r) = 6\delta_{n,t}e^{t/r}.
\end{gather*}
\end{Lem}
\beginpf Take $z \in \D$ such that $|z| = r$. By \eqref{definition of G_n,t} and \eqref{pnt bound}, we have
\begin{align*}
\left|z^{-n}G_{n,t}(z) - P_{n,t}(z)\right| 
&=  \delta_{n,t}|P_{n,t}(z)|\le  \delta_{n,t}\left(|e^{it(z + 1/z)}| + \frac{t^nr^{-n} e^{t/r}}{n!}\right)\\
&\le \delta_{n,t}\left(e^{t(1/r-r)} + \frac{t^nr^{-n} e^{t/r}}{n!}\right)\\
&\le \delta_{n,t}\left(e^{t/r} +\delta_{n,t}r^{-n}e^{t/r}\right) = \delta_{n,t}e^{t/r}(1 +\delta_{n,t}r^{-n}).
\end{align*}
Furthermore, we have
\begin{gather*}
\left|z^{-n}G_{n,t}(z) - e^{it(z + 1/z)}\right|\le \left|z^{-n}G_{n,t}(z) - P_{n,t}(z)\right| + \left|P_{n,t}(z) - e^{it(z + 1/z)}\right|.
\end{gather*}
The last two estimates together with \eqref{p and exp diff} imply
\begin{gather*}
\left|z^{-n}G_{n,t}(z) - e^{it(z + 1/z)}\right| \le \delta_{n,t}e^{t/r}(1 +\delta_{n,t}r^{-n}) + \frac{t^nr^{-n} e^{t/r}}{n!} \le \delta_{n,t}e^{t/r}(1 +2r^{-n}).
\end{gather*}
This gives
$$\max_{|z|=r}|G_{n, t}| \le r^n e^{t(1/r-r)} + 3\delta_{n,t}e^{t/r} \le e^{t/r}(r^n + 3\delta_{n,t}).
$$
So, we have
$$
\max_{|z| = r}\left|zG_{n,t}(z) - z^{n+1}e^{it(z + 1/z)}\right| \le  \delta_{n,t}e^{t/r}r^{n+1}(1 +2r^{-n}) = 
\delta_{n,t}e^{t/r}(r^{n+1} +2r),
$$
and 
$$
\max_{|z| = r}\left|G_{n + 1,t}(z) - z^{n + 1}e^{it(z + 1/z)}\right| \le \delta_{n+1,t}e^{t/r}r^{n+1}(1 +2r^{-(n+1)}) \le \delta_{n,t}e^{t/r}(r^{n+1} +2),
$$
where we used the inequality $\delta_{n+1,t} \le \delta_{n,t}$ for $n>t>0$. It remains to write
\begin{gather*}
\left|G_{n + 1,t}(z) - zG_{n,t}(z)\right| \le \left|G_{n + 1,t}(z) - z^{n + 1}e^{it(z + 1/z)}\right| + \left|zG_{n ,t}(z) - z^{n + 1}e^{it(z + 1/z)}\right|
\end{gather*}
and use the last two estimates. \qed

\medskip

\begin{Lem}\label{l32} For every $n > 0$, $t > 0$, $r \in (0,1)$, we have 
$\sum_{k \ge n} S_k(t,r)r^{-k} \le 6 \delta_{n,t}e^{2t/r} \cdot r^{-n}$.
\end{Lem}
\beginpf For $n > t > 0$ we have $\delta_{n+1, t} = \delta_{n,t}\frac{t}{n + 1}$ hence
\begin{align*}
\sum_{k \ge n} S_k(t,r) r^{-k}
&\le 6\delta_{n,t}e^{t/r}r^{-n}\left(1+\frac{t/r}{n+1} + \frac{(t/r)^2}{(n + 1)(n+2)} + \ldots\right),\\ &\le 6\delta_{n,t}e^{t/r}r^{-n} e^{t/r}  = 6 \delta_{n,t}e^{2t/r} \cdot r^{-n}.
\end{align*}
This is the required estimate. \qed 

\medskip

The following lemma will be proved in Section \ref{s6}, see page \pageref{pfl34}.

\begin{Lem}\label{l34}
Suppose that $q \in \ell^{2}(\Z,\D)$ is such that $\supp q \subset \Z_+$ and let $\fc_{q}$ be defined as in Section \ref{section: Solution for compactly supported initial data}. Then the recurrence coefficients of $\fc_{q}$ coincide with the sequence $\{q(k)\}_{k \ge 0}$.
\end{Lem}

\medskip

\noindent {\bf Proof of Theorem \ref{t3}.} Let $t > 0$, and let $q_0: \Z \to \D$ be a sequence compactly supported on $\Z_+$. Assume that $\prod_{n\in \Z_+}(1-|q_0(n)|^2) \ge \eta$ for some $\eta > 0$. Define the functions $\fc_{q_0} = \bar b/a$, $F_{n,0} = G_{n, t}\fc_{q_0}$ and $F_{n,k}$ as in Section \ref{section: Solution for compactly supported initial data}. Let also $\tilde q_n(t, j) = F_{n,n+j}(0)$, $j\ge -n$, $\tilde q_n(t, j) = 0$, $j< -n$, for $j \in \Z$. We are going to show that $\{\tilde q_n(t, j)\}_{n\ge 0}$ is a Cauchy sequence for each $j \in \Z$. 
Take two positive integers $n_2 > n_1 \ge -j$, fix $r\in(0,1)$ and consider the difference
\begin{align*}
|\tilde q_{n_2}(t, j) - \tilde q_{n_1}(t, j)|
&=|F_{n_2,n_2+j}(0) - F_{n_1,n_1+j}(0)|
\\
&\le \|F_{n_2, n_2+j} - F_{n_1, n_1+j}\|_{L^2(r\T)} 
\le \sum_{k=n_1}^{n_2-1} \|F_{k+1,k+1+j} - F_{k, k+j}\|_{L^2(r\T)}.
\end{align*}
Since $G_{n,t}$ is a contraction by Lemma \ref{lem: modulus gnt}, we have $|G_{k,t}\fc_{q_0}|\le |\fc_{q_0}|$ on $\T$ hence
\begin{gather*}
	\min(\eta(G_{k+1,t}\fc_{q_0}), \eta(zG_{k,t}\fc_{q_0})) \ge \eta(\fc_{q_0}) \ge \eta
\end{gather*}
for every $k$ by Szeg\H{o} theorem \ref{szthm} and our assumption. For a function $F \in \Sch_*(\D)$, denote by $(F)_k$ the $k$-th Schur iterate of $F$ (see \eqref{eq20}, where $(F)_k$ are denoted by $F_k$). Note that $(F)_k = (zF)_{k+1}$. By Theorem~\ref{t0}, we have
\begin{align*}
\|F_{k+1, k+1+j} - F_{k, k+j}\|_{L^2(r\T)}
&=\|(G_{k+1,t}\fc_{q_0})_{k+1+j} - (G_{k,t}\fc_{q_0})_{k+j}\|_{L^2(r\T)}\\
&=\|(G_{k+1,t}\fc_{q_0})_{k+1+j} - (zG_{k,t}\fc_{q_0})_{k+1+j}\|_{L^2(r\T)}\\
&\le 
C(\eta, r)r^{-k-1-j}\|G_{k+1,t}\fc_{q_0} - zG_{k,t}\fc_{q_0}\|_{L^2(r\T)}
\\
&\le C(\eta, r)r^{-k-1-j}\|G_{k+1,t} - zG_{k,t}\|_{L^2(r\T)}.
\end{align*}
Using Lemma \ref{Lemma: convergence of gnt} for $n_1  > t >0$ such that $\delta_{n_1,t} < 1$, we can proceed as follows:
\begin{align*}
\|F_{k+1, k+1+j} - F_{k, k+j}\|_{L^2(r\T)}
\le C(\eta, r)r^{-k-j-1}\max_{|z| = r}|G_{k+1,t} - zG_{k,t}| \le C(\eta, r)S_k(t,r)r^{-k-j-1}.
\end{align*}
From Lemma \ref{l32} we now see that
\begin{align*}
|\tilde q_{n_2}(t, j) - \tilde q_{n_1}(t, j)|
&\le 
r^{-j-1} C(\eta, r) \cdot \sum_{k=n_1}^{\infty}S_k(t,r)r^{-k}\le 6C(\eta, r) \delta_{n_1,t}e^{2t/r} \cdot r^{-n_1-j-1}.
\end{align*}
Recall that $\delta_{n_1,t} = \frac{e^tt^{n_1}}{n_1!}$ decays very rapidly as $n_1\to\infty$, thus, $\{\tilde q_{n}(t, j)\}_{n\ge -j}$ is a Cauchy sequence for every $j \in \Z$. Denote its limit by $\tilde q(t, \cdot)$. 
Letting $n_1 = n$ and taking the limit in as $n_2 \to +\infty$, we obtain  
\begin{gather*}
|\tilde q(t, j) - \tilde q_n(t, j)| 
\le 6C(\eta, r) \delta_{n,t}e^{2t/r} \cdot r^{-n-j-1}.
\end{gather*}
Taking $r = 1/2$ (any other $r\in (0,1)$ will do) and using the inequality $n! \ge \sqrt{2\pi n}(n/e)^n$, we get
\begin{align*}
|\tilde q(t, j) - \tilde q_n(t, j)| \le 6C(\eta, 1/2)\frac{e^t t^n}{n!}e^{4t}2^{n + j + 1}
= 2^{j}C(\eta, 1/2)\frac{12e^{5t}}{\sqrt{2\pi n}}\left(\frac{2et}{n}\right)^{n},
\end{align*}
where $n \in \Z_+$, $j \in \Z$, $t>0$ are such that $n+j \ge 0$, $n > t > 0$, and $\delta_{n,t} <1$. 

It remains to show that 
$\tilde q(t, j) = q(t, j)$, i.e., $\tilde q$ solves Ablowitz-Ladik equation \eqref{al} with the initial datum $q_0$. By assertions $(6)$, $(7)$ of Theorem \ref{t5} it is suffices to check that $\rc_{\tilde q} = \rc_{q}$, equivalently, $\rc_{\tilde q} = e^{it(z+1/z)}\rc_{q_0}$.

Note that $\tilde q_n(t, \cdot - n)$ is supported on $\Z_+$, moreover, we have 
$\tilde q_n(t, j - n) = F_{n, j}(0)$ for $j \in \Z_+$. Let us denote the coefficients  in \eqref{eq25bisbis} for $q_{0}$, $\tilde q_n(t, \cdot - n)$, by $a$, $b$, and $a_{n, 0}$, $b_{n,0}$, respectively. We have $\rc_{\tilde q_n(t, \cdot - n)} = b_{n,0}/a_{n, 0}$ and $F_{n, 0} = \fc_{\tilde q_n(t, \cdot - n)} = \ov{b_{n, 0}}/a_{n, 0}$, where equality $F_{n, 0} = \fc_{\tilde q_n(t, \cdot - n)}$ holds by Lemma \ref{l34} because the recurrence coefficients of $F_{n,0}$ coincide with the sequence $\{\tilde q(t,j-n)\}_{j \ge 0}$. By \eqref{p and exp diff}, we also have 
\begin{gather*}
	z^{-n}\ov{b_{n,0}}/a_{n,0} = z^{-n}F_{n,0} = z^{-n}G_{n,t}\fc_{q_0} \to e^{it(z+1/z)}\fc_{q_0} = e^{it(z+1/z)}\bar b/a
\end{gather*}
uniformly on $\T$. We now will use well-known properties of coefficients $a$, $b$ in \eqref{eq25bisbis}.
 Namely, the functions $a$, $a_{n, 0}$ are outer, have positive values at $z = 0$, and satisfy $|a|^2 - |b|^2 = 1$, $|a_{n, 0}|^2 - |b_{n, 0}|^2 = 1$ on $\T$ (for the proof, see Section~\ref{s6}). Convergence $z^{-n}\ov{b_{n,0}}/a_{n,0} \to e^{it(z+1/z)}\bar b/a$ then implies $|a_{n, 0}|^2 \to |a|^2$, $\log|a_{n, 0}|^2 \to \log|a|^2$ uniformly on $\T$, hence $a_{n, 0} \to a$ in Lebesgue measure on $\T$ by properties of outer functions (more precisely, by the weak continuity of the Hilbert transform, see a discussion next to formula \eqref{eq11new}). It follows that $z^n b_{n, 0} \to e^{it(z+1/z)} b$ in Lebesgue measure on $\T$. Therefore,
\begin{equation}\label{eq40}
\rc_{\tilde q_n(t, \cdot)} = z^{n}\rc_{\tilde q_n(t, \cdot - n)} = z^{n}b_{n,0}/a_{n,0} \to  e^{it(z+1/z)} b/a = e^{it(z+1/z)} \rc_{q(t, \cdot)},
\end{equation}
in Lebesgue measure on $\T$ (the first equality in \eqref{eq40} is assertion $(3)$ of Theorem \ref{t5}). 
On the other hand, as $n \to +\infty$ the quantities
$$
\int_{\T}\log(1-|\rc_{\tilde q_n(t, \cdot)}|^2)\,dm = \int_{\T}\log(|a_{n,0}|^{-2})\,dm = \log|a_{n,0}(0)|^{-2}
$$
tend to
$$
\log|a(0)|^{-2} =\int_{\T}\log(1-|\rc_{q}|^2)\,dm = \int_{\T}\log(1-|e^{it(z+1/z)}\rc_{q}|^2)\,dm.
$$
Then, taking into account \eqref{eq40}, we see that $\rc_{\tilde q_n(t, \cdot)} \to e^{it(z+1/z)}\rc_{q}$ in the metric space $X$ by Proposition \ref{prop: convergence in measure to convergence in X}. Moreover, the quantities $\esssup_{\T}(|\rc_{\tilde q_n(t, \cdot)}|^2) = 1-\esssup_{\T}|a_{n, 0}|^{-2}$ are uniformly separated from $1$ because $a_{n, 0}$ converge uniformly on $\T$ to the bounded function $a$. Then continuity of the inverse NLFT map (i.e., assertion $(6)$ in Theorem \ref{t5}) gives us the convergence of $\tilde q_n(t, \cdot)$ to $\F_{sc}^{-1}(e^{it(z+1/z)}\rc_{q})$ in a subspace $\cG(\delta)$, $\delta \in (0,1)$, of the metric space $\ell^2(\Z, \D)$. Since the sequence $\tilde q_n(t, \cdot)$ converges elementwise to $\tilde q(t, \cdot)$ as $n \to +\infty$, we get $\tilde q(t, \cdot) = \F_{sc}^{-1}(e^{it(z+1/z)}\rc_{q})$ on $\Z$. Then, $\rc_{\tilde q(t, \cdot)} = \F_{sc}(\tilde q(t, \cdot)) = e^{it(z+1/z)}\rc_{q}$ 
almost everywhere on $\T$, and the proof is completed. \qed 

\begin{Rema}\label{r31} In the proof of Theorem \ref{t3}, we have used the fact that \eqref{al} is solvable for compactly supported initial data. This can be proved by a variety of methods, see Appendix for a direct proof in a much more general situation. Assertions $(6)$, $(7)$ in Theorem \ref{t5} guarantee that the solution will be determined by its reflection coefficient $\rc_{q(t, \cdot)} = e^{it(z+1/z)}\rc_{q_0}$ at any moment of time $t \in \R$. 
\end{Rema}

\section{Localization. Proof of Theorem \ref{t2}}
The following lemma is well-known, see, e.g., (1.3.43) in \cite{Simonbook1}. 
\begin{Lem}\label{lemma: schur functions witl the same an}
Let $F,G\in \Sch_*(\D)$, and let $F_k$, $G_k$ be their Schur iterates \eqref{eq20}. Assume that $F_k(0) = G_k(0)$ for $0 \le k \le n$. Then
$\max_{|z| =r}|F(z) - G(z)|\le 2r^{n+1}$.
\end{Lem}
\begin{Lem}
\label{from diff of Schur functions to diff of solutions}
Let $F, G \in \Sch(\D)$ be such that $\min(\eta(F),\eta(G)) \ge \eta$ for some $\eta > 0$.  Denote by $F_k$, $G_k$ their Schur iterates \eqref{eq20}, and consider the solutions of \eqref{al} with the initial value 
\begin{gather*}
q_{0,F} = 
\begin{cases}
	F_n(0), &n\ge 0,\\
	0, &n < 0,
\end{cases}
\qquad
q_{0,G} = 
\begin{cases}
	G_n(0), &n\ge 0,\\
	0, &n < 0.
\end{cases}
\end{gather*}
Denote them by $q_F$ and $q_G$, respectively.
Then for every $n>t>0$, $r\in (0,1)$, the inequality
\begin{gather*}
|q_F(t, j) - q_G(t, j)| \le r^{-j} e^{t/r} C(\eta, r)\sup_{|z| = r}|F(z) - G(z)|, 
\end{gather*}
holds for all $j\in \Z$. Here $C(\eta, r)$ is the function from Theorem \ref{t0}.
\end{Lem}
\beginpf For a function $H \in \Sch_*(\D)$, let us denote by $(H)_k$ its Schur iterates \eqref{eq20}. By Theorem \ref{t3}, we have
\begin{gather*}
q_F(t, j) = \lim_{n\to\infty} (G_{n,t}F)_{n + j}(0), \qquad q_G(t, j) = \lim_{n\to\infty} (G_{n,t} G)_{n + j}(0),\qquad j\in \Z.
\end{gather*}
Therefore, we can apply Theorem \ref{t0} and the bound $|G_{n,t}| < e^{t/r}(r^n + 3\delta_{n,t})$ from Lemma \ref{Lemma: convergence of gnt} to get
\begin{align}
|q_F(t, j) - q_G(t, j)|&\le \limsup_{n\to\infty} |(G_n F)_{n + j}(0) - (G_n G)_{n + j}(0)|\notag
\\
&\le \limsup_{n\to\infty} \|(G_{n,t}F)_{n + j} - (G_{n,t}G)_{n + j}\|_{L^2(r\T)} \notag \\
&\le \limsup_{n\to\infty} C(\eta, r)r^{-n - j}\|G_{n,t}F - G_{n,t}G\|_{L^2(r\T)}
\\
&\le \limsup_{n\to\infty} C(\eta, r)r^{-j}e^{t/r}(1 + 3\delta_{n,t}r^{-n})\sup_{|z| = r}|F(z) - G(z)|\label{eq42}\\
&=r^{-j} e^{t/r} C(\eta, r)\sup_{|z| = r}|F(z) - G(z)|,\notag
\end{align}
where we have used in \eqref{eq42} the convergence $\delta_{n,t}r^{-n} \to 0$ as $n \to \infty$. \qed

\medskip

\noindent {\bf Proof of Theorem \ref{t2}.} Recall that $q_0:\Z \to \D$ is such that $\prod_{n\in \Z}(1 - |q_0(n)|^2)\ge \eta > 0$, the sequence $q_{0, N}$ is defined by
\begin{gather*}
q_{0, N}(n) = 
\begin{cases}
q_{0}(n), &|n|\le N,\\
0, & |n|> N,
\end{cases}
\end{gather*}
and $q_{N}$ is the corresponding solution of \eqref{al} (see Remark \ref{r31}). 
Let $C(\eta, r)$ be the function from Theorem \ref{t0}. We want to prove the inequality
\begin{gather}\label{final convergence estimate }
|q_{N+K}(t, j) - q_N(t,j)|\le \frac{4e^{t/r} C(\eta, r)}{1 - r} r^{N - |j|}, \qquad K \in \Z_+.
\end{gather}
Then $\{q_{N}(t, j)\}$ will be a Cauchy sequence for each $t, j$, and its limit, to be denoted by $q$, solves \eqref{al}. This is easy to check if one rewrites \eqref{al} in the integral form. Estimate \eqref{eq: theorem 1.3 ineq} will follow from \eqref{final convergence estimate } by taking the limit as $K \to +\infty$.

\smallskip

\noindent For integer numbers $A\le B$, consider the sequences $q_{0,[A,B]}$, $\refl {q}_{0,[A,B]}$ in $\ell^2(\Z, \D)$ defined by
\begin{gather*}
q_{0,[A,B]}(j) = q_0(j + A)\mathbf{1}_{[0, B - A]}(j), \qquad \refl q_{0,[A,B]}(j) = q_0(-j + B)\mathbf{1}_{[0, B - A]}(j),
\end{gather*}
where $\mathbf{1}_{S}$ is the indicator function of a set $S$. These sequences both supported on $[0, B - A]$ and their entries are symmetric on this segment. Denote the corresponding solutions of \eqref{al} by $q_{[A,B]}$ and $\refl {q}_{[A,B]}$. By properties of \eqref{al}, the symmetry relation
\begin{gather}\label{eq: symmetry for AB}
q_{[A, B]}(t, j) = \refl q_{[A, B]} (t, B - A - j), \quad t\in \R, \quad j\in \Z,
\end{gather}
holds for each $t \in \R$. Moreover, comparing this with the definition of $q_N$, we see that $q_N(t, j) = q_{[-N, N]}(t, j + N).$ The inequality \eqref{final convergence estimate } will follow by summing up a telescoping series if we check the estimate
\begin{gather}
\label{eq: cauchy conv for q[N]}
|q_N(t,j) - q_{N + 1}(t,j)|\le 4C(\eta, r)e^{t/r}r^{N-|j|}, \qquad N\ge |j|. 
\end{gather}
In the new notation the latter takes the form
\begin{gather}\label{rewritten final estimate}
|q_{[-N, N]}(t, j + N) - q_{[-N - 1, N + 1]}(t, j + N + 1)|\le 4C(\eta, r)e^{t/r}r^{N-|j|}, \quad N\ge |j|. 
\end{gather}
For $A \le B$ let $f_{[A, B]}$ and $\refl f_{[A, B]}$ be the Schur functions which recurrence coefficients are $q_{0,[A,B]}\mid_{\Z_+}$ and $\refl {q}_{0,[A,B]}\mid_{\Z_+}$ respectively. The Schur functions $f_{[-N, N]}$, $f_{[-N, N + 1]}$ have the same first $2N + 1$ Schur coefficients. Hence by Lemmas~\ref{lemma: schur functions witl the same an} and \ref{from diff of Schur functions to diff of solutions} we get
\begin{align}
\notag
\left|q_{[-N, N]}(t,n) - q_{[-N, N + 1]}(t,n)\right|&\le  r^{-n} \cdot C(\eta, r)e^{t/r}\sup_{|z| = r}\left|f_{[-N, N]}(z) - f_{[-N, N + 1]}(z)\right|
\\
\label{fisrt part of the final}
&\le 2C(\eta, r)e^{t/r} r^{2N-n + 1}, 
\end{align}
for all $n \ge 0$. Similarly, the functions $\refl f_{[-N, N + 1]}$ and $\refl f_{[-N - 1, N + 1]}$ have coinciding first $2N+1$ Schur coefficients, therefore
\begin{gather}\label{second part of the final}
\left|\refl q_{[-N, N + 1]}(t,n) - \refl q_{[-N - 1, N + 1]}(t,n)\right|\le 2C(\eta, r)e^{t/r} r^{2N-n + 1},\qquad n\ge 0.
\end{gather}
Notice that
\begin{multline*}
|q_{[-N, N]}(t, n) - q_{[-N - 1, N + 1]}(t, n + 1)| \le 
\\
\le |q_{[-N, N]}(t, n) - q_{[-N, N + 1]}(t, n)| + |q_{[-N, N + 1]}(t, n) - q_{[-N - 1, N + 1]}(t, n + 1)| \le \\
\le 2C(\eta, r)e^{t/r} r^{2N-n + 1} + |q_{[-N, N + 1]}(t, n) - q_{[-N - 1, N + 1]}(t, n + 1)|.
\end{multline*}
By relation \eqref{eq: symmetry for AB}, the last term equals
\begin{align*}
|\refl q_{[-N, N + 1]}(t, 2N + 1 - n) - \refl q_{[-N - 1, N + 1]}(t, 2N + 1- n)| 
&\le 2C(\eta, r)e^{t/r} r^{2N-(2N + 1 - n) + 1}\\
&= 2C(\eta, r)e^{t/r} r^{n},
\end{align*}
where we used \eqref{second part of the final} in the first inequality. 
Therefore, we have $|q_{[-N, N]}(t, n) - q_{[-N - 1, N + 1]}(t, n + 1)| \le 2C(\eta, r)e^{t/r} (r^{2N-n + 1} + r^{n}).$ 
Substitution of $n=j + N$ then gives 
\begin{align*}
|q_{[-N, N]}(t, j + N) - q_{[-N - 1, N + 1]}(t, j + N + 1)|
&\le 2C(\eta, r)e^{t/r} (r^{N-j + 1} + r^{j + N})\\
&\le 4C(\eta, r)e^{t/r}r^{N-|j|},
\end{align*}
which is \eqref{rewritten final estimate}. 
\qed

\section{Complexity of the algorithm}\label{s7}
In the introduction, we claimed that the algorithm outlined in Section \ref{s1-4} takes $O(\mathbf{n}\log^2\mathbf{n})$ operations for $\mathbf{n} = t+\log\eps^{-1}$. Here we prove this estimate. 

\medskip

Let $q_0 \in \ell^2(\Z, \D)$ be such that $\prod_{n \in \Z}(1-|q_0(n)|^2) \ge \eta > 0$, and let $t>0$. Take $\eps \in (0,1)$, set $r = 1/2$, and choose $N \in \Z_+$ such that  
the right hand side in \eqref{eq: theorem 1.3 ineq} does not exceed $\eps/2$ at $j=0$:
$$
8e^{2t}C(\eta, 1/2) 2^{-N} \le \eps/2, \qquad C(\eta, 1/2) = \exp\left(17 \log\eta^{-1}\cdot \left(2 + \frac{1}{1-\sqrt{1-\eta}}\right)\right).
$$
Since $8e^{2t}C(\eta, 1/2) 2^{-N} \le 2^{-N+4+3t}C(\eta, 1/2)/2$, one can take any $N \ge 5 + [3t + \log_2 \frac{C(\eta, 1/2)}{\eps}]$.
Then, choose the window $\Delta = [n_0 - N, n_0 + N]$, truncate $q_0$ by setting $q_0 = 0$ on $\Z \setminus \Delta$, and shift $q_0$ by $n_0 - N$ to make it supported in $[0, 2N]$. Denote the resulting sequence by $q_{0, [n_0 - N, n_0 + N]}$. Choose $n > t$ so that $\delta_{n,t} <1$ and
$$
2^{j}C(\eta, 1/2)\frac{12e^{5t}}{\sqrt{2\pi n}}\left(\frac{2et}{n}\right)^{n} \le \eps/2, \quad j = N.
$$
Since we already have $8e^{2t}C(\eta, 1/2) 2^{-N} \le \eps/2$, it suffices to choose $n$ so that 
$$
2^{2N}e^{3t} \frac{12}{8\sqrt{2\pi n}}\left(\frac{2et}{n}\right)^{n} \le 1.
$$
For $n \ge 2N \ge 8et \ge 5t$, we have
$$
2^{2N}e^{3t} \frac{12}{8\sqrt{2\pi n}}\left(\frac{2et}{n}\right)^{n} \le 2^{2N + 5t}\left(\frac{2et}{n}\right)^{n} \le 2^{5t}\left(\frac{2et}{N}\right)^{2N} \le \left(\frac{4et}{N}\right)^{2N} \le 1,
$$  
therefore, one can take $n = 2N$, $N = 5 + [4et + \log_2 \frac{C(\eta, 1/2)}{\eps}]$. Note that with this choice 
$$
\delta_{n, t} 
= \frac{t^n e^t}{n!} \le \left(\frac{te}{n}\right)^n e^t \le \left(\frac{te^2}{n}\right)^n \le \left(\frac{8et}{2N}\right)^n < 1.
$$
We see that for $n = 2N$, $N = 5 + [4et + \log_2 \frac{C(\eta, 1/2)}{\eps}]$, Theorem \ref{t3} applied to $q_{0, [n_0 - N, n_0 + N]}$ in place of $q_0$ will give a sequence $\tilde q_n$ approximating the corresponding solution $q_{[n_0 - N, n_0 + N]}$ with accuracy $|\tilde q_n(t,N) - q_{[n_0 - N, n_0 + N]}(t, N)| \le \eps/2$. Then $|\tilde q_n(t,N+1) - q(t, n_0)| < \eps$ and it remains to estimate the number of operations that are needed to construct $\tilde q_n(t,N)$ from $q_0$ for $n = 2N$. 

\smallskip

Having $q_0$, $t_0$, $n_0$, $\eps$, $\eta$, we set $N = 5 + [4et + \log_2 \frac{C(\eta, 1/2)}{\eps}]$ and define array $q_{0, [n_0 - N, n_0 + N]}$ of $2N+1$ elements. Then we use formula \eqref{eq25bisbis} to find $a$, $b$. This can be done either by a direct multiplication of $2N+1$ matrices in $O(N^2)$ operations or by using a dyadic divide-and-conquer multiplication algorithm together with the fast Fourier transform (FFT) in $O(N\log^2 N)$ operations. Next, define coefficients of polynomials $P = G_{n,t}\bar b$, $Q = a$ (two arrays of length $2n+1+2N+1$, $2N+1$, respectively). This takes $O(N^2)$ operations in naive realization of multiplications of polynomials or $O(N\log N)$ operations with FFT. Taking $n+N+1$ steps of Schur's algorithm for $P/Q$, we find $\tilde q_{n}(t, j)$ on $[0,N]$, which solves the problem. Straightforward  realization of Schur's algorithm based on its definition requires $O(N^2) = O(\mathbf{n}^2)$ operations (recall that $\mathbf{n} = \log\eps^{-1} + t$). It could be fastened up to $O(\mathbf{n}\log^2 \mathbf{n})$ operations with more delicate realization, see Section 2.2 in \cite{Ammar89}. Notice that the numerical experiments in \cite{Ammar89} use arithmetic of real numbers, while complexity estimate $O(\mathbf{n} \log^2 \mathbf{n})$ given on page 192 in \cite{Ammar89} holds for complex data. As reader can see from the algorithm, the same $O(\mathbf{n} \log^2 \mathbf{n})$ operations (with worsted constant) are sufficient to find $\tilde q_{n}(t,N)$ on $[0,2N]$ and approximate $q(t,\cdot )$ with accuracy $\eps$ on the interval $[n_0-N/2, n_0+N/2]$, not only at the point $n_0$. It is also worth mentioning that the question of numerical stability (in our case -- estimating round-off errors and taking into account issues related to arithmetic of long numbers) deserves a special consideration, it does not treated neither in \cite{Ammar89} nor in this paper.       

\medskip

\section{The nonlinear Fourier transform. Proof of Theorem \ref{t5}}\label{s6}
In this section we collect some basic facts about the nonlinear Fourier transform (NLFT). Some of them were used in the first part of the paper. The reader can find more information in the preprint \cite{TT} or in its extended version \cite{TTT}. 

\medskip

The exposition in this section is independent from the first part of the paper. Let us recall the definition of the NLFT map for the reader's convenience. For $p\ge 1$, define $\ell^p(\Z, \D)$ as a set of sequences $q\colon \Z\to\D$ satisfying $|q(n)| < 1$ for every $n \in \Z$ and $\sum_{n \in \Z} |q(n)|^p < \infty$. The set $\ell^p(\Z_+, \D)$ is defined similarly with $\Z_+ = \Z \cap [0, +\infty)$. Take a sequence $q \in \ell^1(\Z, \D)$ and define $a$, $b$ by
\begin{equation}\label{eq25bis}
\begin{pmatrix}
a & b\\
\ov{b} & \bar a
\end{pmatrix} 
=
\prod_{k \in \Z}
\frac{1}{\sqrt{1-|q(k)|^2}}
\cdot
\begin{pmatrix}
1 & \ov{q(k)} z^{-k}\\
q(k)  z^k & 1
\end{pmatrix}, \qquad z \in \T.
\end{equation}  
Here, the product $\prod_{k\in \Z}T_k$ of matrices $T_k$ is understood as the limit 
$\lim_{n \to +\infty}T_{-n}T_{-n+1}\cdot\ldots\cdot T_{n-1}T_{n}$. Assumption $q \in \ell^1(\Z, \D)$ guarantees that the product converges uniformly on $\T$. 
We will see in Section~\ref{s52} that the product in \eqref{eq25bis} has the form $\left(\begin{smallmatrix}
a & b\\
\ov{b} & \bar a
\end{smallmatrix}\right)$ for some $a$, $b$. The authors of \cite{TT} define NLFT as the map that sends $q$ to the pair $\left(\begin{matrix}
a & b\end{matrix}\right)$. We will use an equivalent definition and consider the so-called reflection coefficient $\rc_{q} = \frac{b}{a}$ in place of $\left(\begin{matrix}
a & b\end{matrix}\right)$. So, in our case, NLFT takes $q$ into $\rc_q$. In the next two subsections we define the reflection coefficient as an object of the theory of orthogonal polynomials on the unit circle. We also prove equivalence of the two definitions of NLFT map.  
\subsection{Szeg\H{o} measures and Szeg\H{o} functions} Let $\mu$ be a probability measure supported on an infinite subset of the unit circle $\T = \{z \in \C: |z|=1\}$ of the complex plane, $\C$. For $n \in \Z_+$, denote by $\Phi_n$ the monic orthogonal polynomial of degree $n$ generated by $\mu$, and set $\Phi^*_n = z^n \ov{\Phi_n(1/\bar z)}$. These polynomials satisfy the following relation:
\begin{gather}\label{recurrence relation}
	\Phi_{n+1} = z\Phi_n - \bar \alpha_n \Phi_n^*, \quad n \ge 0, \quad \Phi_0 = 1,
\end{gather}
where the {\it recurrence coefficients}, $\alpha_n$, $n \ge 0$, lie in the open unit disk $\D = \{z \in \C: \; |z|<1\}$. Conversely, any sequence $\{\alpha_n\}_{n \ge 0} \subset \D$ gives rise to a unique probability measure 
$\mu$ on $\T$ whose closed support $\supp \mu$ contains infinitely many points. These two facts can be found in Section 1.7 of \cite{Simonbook1}. The Schur function $f$ of a probability measure $\mu$ on $\T$ is defined by
\begin{equation}\label{eq0}
\frac{1 + zf(z)}{1-zf(z)} = \int_{\T}\frac{1+\bar \xi z}{1-\bar \xi z}\,d\mu(\xi), \qquad z \in \D.
\end{equation}
Notice that \eqref{eq0} provides a bijective correspondence between Schur functions and measures on $\T$. Taking the real part in both sides of this equality, we get
\begin{equation}\label{eq4}
\frac{1 - |zf(z)|^2}{|1-zf(z)|^2} = \int_{\T}\frac{1-|z|^2}{|1-\bar \xi z|^2}\,d\mu(\xi), \qquad z \in \D. 
\end{equation}
From \eqref{eq0}, \eqref{eq4}, and Schwarz lemma we see that $f$ indeed belongs to the Schur class $\Sch(\D)$, i.e., it is analytic in $\D$ and satisfies $\sup_{z \in \D}|f(z)| \le 1$. Recall that the Schur iterates of $f = f_0$ are defined by
\begin{gather}\label{eq: NFLT schur algorithm}
	zf_{n+1} = \frac{f_n - f_n(0)}{1 - \ov{f_n(0)}f_n}, \qquad n \ge 0.
\end{gather}
Geronimus theorem says that recurrence coefficients $\alpha_n$ in \eqref{recurrence relation} coincide with recurrence coefficients in Schur's algorithm:  $\alpha_n = f_n(0)$, $n \ge 0$. See Chapter 3 in \cite{Simonbook1} for the proof. 

\medskip

Let $\mu = w\,dm+\mus$ be the Radon-Nikodym decomposition of $\mu$ into the absolutely continuous and singular parts, where $m$ is the Lebesgue measure on $\T$ normalized by $m(\T) = 1$. Denote by $\{\alpha_n\}$ the set of recurrence coefficients of the measure $\mu$ and let $f$ be its Schur function. An extended version of Szeg\H{o} theorem (Theorem \ref{szthm}) says that conditions $\log w \in L^1(\T)$, $\log(1-|f|^2) \in L^1(\T)$,  $\{\alpha_n\} \in \ell^2(\Z_+,\D)$ are equivalent, and, moreover, 
\begin{equation}\label{eq3}
\int_{\T}\log w(\xi)\,dm(\xi) = \int_{\T}\log(1 - |f(\xi)|^2)\,dm(\xi) = \log\prod_{n \ge 0}(1 - |\alpha_n|^2).
\end{equation}
It is not difficult to see that the three quantities in \eqref{eq3} are defined for any triple $\mu$, $f$, $\{\alpha_n\}$, but could be $-\infty$. In fact, Szeg\H{o} theorem implies that quantities in \eqref{eq3} are finite (i.e., $> -\infty$) or not simultaneously. Measures of Szeg\H{o} class
$$
\szc = \bigl\{\mu = w\,dm + \mus: \; \mu(\T) = 1, \; \log w \in L^1(\T)\bigr\}
$$ 
and their orthogonal polynomials have many interesting properties that constitute rich Szeg\H{o} theory. We will use its part related to a discrete scattering. For this we will need the notion of the dual orthogonality measure, the Szeg\H{o} function, and the dual Szeg\H{o} function. 

\medskip

Consider a probability measure $\mu$ on $\T$ with infinite support. Let, as before, $f$ denote the Schur function of $\mu$. The dual measure $\mu_d$ is defined as the probability measure on $\T$ corresponding to the Schur function $-f$:
\begin{equation}\label{eq0d}
\int_{\T}\frac{1+\bar \xi z}{1-\bar \xi z}\,d\mu_d(\xi) = \frac{1 + zf_d(z)}{1-zf_d(z)} = \frac{1 - zf(z)}{1+zf(z)}, \qquad z \in \D.
\end{equation}
It is not difficult to check that if $\{\alpha_n\}_{n \ge 0}$ is the sequence of recurrence coefficients of $\mu$, then  $\{-\alpha_n\}_{n \ge 0}$ is the sequence of recurrence coefficients of $\mu_{d}$. 
Monic orthogonal polynomials for $\mu_d$ will be denoted by $\Psi_n$. We also will need the normalized orthogonal polynomials for $\mu$ and $\mu_{d}$:
\begin{equation}\label{eq31}
\phi_{n} = \frac{\Phi_{n}}{\|\Phi_n\|_{L^2(\mu)}}, \qquad \phi^*_{n} = \frac{\Phi^*_{n}}{\|\Phi^*_n\|_{L^2(\mu)}}, \quad 
\psi_{n} = \frac{\Psi_{n}}{\|\Psi_n\|_{L^2(\mu_d)}}, \qquad \psi^*_{n} = \frac{\Psi^*_{n}}{\|\Psi^*_n\|_{L^2(\mu_d)}}.
\end{equation}
In fact, 
\begin{equation}\label{eq43}
\|\Phi_n\|_{L^2(\mu)}^2 = \|\Phi^*_n\|_{L^2(\mu)}^2 = 
\|\Psi_n\|_{L^2(\mu_d)}^2 = \|\Psi^*_n\|_{L^2(\mu_d)}^2 = \prod_{k=0}^{n-1}(1-|\alpha_k|^2), 
\end{equation}
for all $n \ge 1$, see Chapter 3.2 in  \cite{Simonbook1}. The Szeg\H{o} function, $D_{\mu}$, of a measure $\mu = w\,dm + \mus$ from Szeg\H{o} class $\szc$ is the outer function in the open unit disk $\D$ such that $D_{\mu}(0)>0$ and $|D_{\mu}|^{2} = w$ Lebesgue almost everywhere on $\T$ in the sense of nontangential boundary values. It could be defined by the formula
\begin{gather}\label{D in z formula}
	D_{\mu}(z) = \exp\left(\frac{1}{2}\int_{\T}\frac{1+\bar \xi z}{1-\bar \xi z}\log w(\xi)\, dm(\xi)\right), \qquad z \in \D.
\end{gather}
It follows from the Szeg\H{o} theorem (see \eqref{eq3}) that $\mu \in \szc$ if and only if $\mu_d \in \szc$. We will denote the Szeg\H{o} function of $\mu_d$ by $D_{\mu_d}$. 
It is known that $\phi_n^*\to  D_{\mu}^{-1}$, $\psi_n^*\to D_{\mu_d}^{-1}$ as $n\to\infty$ in $\D$ and
\begin{gather}\label{eq: ratio of Szego functions equals schur }
	\frac{1 + zf}{1 - zf} = \lim_{n \to \infty}\frac{\Psi_n^*(z)}{\Phi_n^*(z)} =\lim_{n \to \infty}\frac{\psi_n^*(z)}{\phi_n^*(z)} =\frac{D_{\mu_d}^{-1}(z)}{D_{\mu}^{-1}(z)}, \qquad z \in \D,
\end{gather}
see Theorem 2.4.1 and Chapter 3.2 in  \cite{Simonbook1}. In particular, we have
\begin{equation}\label{eq5}
\Re\left(D_{\mu_d}^{-1}\ov{D_{\mu}^{-1}}\right) = \Re\left(\frac{D_{\mu_d}^{-1}}{D_{\mu}^{-1}}\right)|D_{\mu}|^{-2} = \frac{1 - |z f|^2}{|1 - z f|^2}|D_{\mu}|^{-2} = w|D_{\mu}|^{-2} = 1
\end{equation}
almost everywhere on $\T$ in the sense of non-tangential boundary values.

\medskip

\subsection{Reflection coefficients}\label{s52} Let us now define a reflection coefficient of a sequence $q$ in $\ell^2(\Z, \D)$. To simplify notation, we set $q_n = q(n)$, $n \in \Z$. Consider the sequences $\{\alpha_n\}_{n \in \Z_+}$ and $\{\beta_n\}_{n \in \Z_+}$ from $\ell^2(\Z_+, \D)$ defined by $\alpha_n = q_n$, for $n\ge 0$ and $\beta_0 = 0$, $\beta_n = -\ov{q_{-n}}$ for $n\ge 1$,
\begin{gather}\label{eq48}
\begin{matrix}
& & & & \alpha_0 & \alpha_1 & \alpha_2 & \alpha_3 & \ldots
\\
\ldots & q_{-3} & q_{-2} & q_{-1} & q_{0} & q_1 & q_2 & q_3& \ldots
\\
\ldots & \beta_3 & \beta_2 & \beta_1 & 0. & & &
\end{matrix}
\end{gather}
Define the measures $\mu^{+}$, $\mu^-$ with the recurrence coefficients $\{\alpha_n\}_{n \ge 0}$, $\{\beta_n\}_{n \ge 0}$, respectively. Let also $\mu_{d}^{\pm}$ be the dual measures corresponding to $\mu^{\pm}$.
Define the Wall analytic functions in $\D$ by
\begin{equation}\label{eq49}
\fa^{\pm} = \frac{D_{\mu_d^{\pm}}^{-1} + D_{\mu^{\pm}}^{-1}}{2},
\qquad 
\fb^{\pm} = \frac{D_{\mu_d^{\pm}}^{-1}-D_{\mu^{\pm}}^{-1}}{2z}.
\end{equation}
The fact that $D^{-1}_{\mu_{d}^{\pm}}(0) = D^{-1}_{\mu^{\pm}}(0)$ follows from \eqref{eq: ratio of Szego functions equals schur }. 
Using \eqref{eq5}, we obtain $|\fa^{\pm}|^2 - |\fb^{\pm}|^2 = \Re\left(D_{\mu_d^\pm}^{-1}\ov{D_{\mu^\pm}^{-1}}\right) = 1$ Lebesgue almost everywhere on $\T$ in the sense of non-tangential boundary values. Also, we have 
\begin{equation}\label{eq44}
\frac{1+z\frac{\fb^{\pm}}{\fa^{\pm}}}{1-z\frac{\fb^{\pm}}{\fa^{\pm}}} = \frac{D_{\mu_d^{\pm}}^{-1}(z)}{D_{\mu^{\pm}}^{-1}(z)} = \frac{1 + zf^{\pm}}{1 - zf^{\pm}},
\end{equation} 
for the Schur functions $f^{\pm}$ of $\mu^{\pm}$, hence $f^{\pm} = \fb^\pm/\fa^\pm$. On $\T$, we set
\begin{equation}\label{eq26}
a = \fa^+\fa^{-} - \fb^+\fb^-, \qquad b = \fa^-\ov{\fb^+} - \fb^-\ov{\fa^+}.
\end{equation}
Below we will use the fact that $a$ is defined by \eqref{eq26} not only on $\T$ but also in $\D$ and is analytc there. Note that 
$
|a|^2 - |b|^2 = 
(|\fa^{+}|^2 - |\fb^{+}|^2)(|\fa^{-}|^2 - |\fb^{-}|^2)=1
$ 
almost everywhere on $\T$. Next, define the reflection coefficient, $\rc_{q}$, of the sequence $q = \{\alpha_n\}_{n \in \Z}$ in $\ell^2(\Z, \D)$ by
\begin{equation}\label{eq32} 
\rc_{q} = \frac{b}{a}.
\end{equation}
It is possible to associate with $q$ an operator on $\ell^2(\Z) \oplus \ell^2(\Z)$  in a way that will place the reflection coefficient $\rc_q$ into the setting of a discrete scattering theory, see \cite{TTT}. Our first proposition collects the properties of objects defined in the present section.
\begin{Prop}\label{pr2}
	For every $q \in \ell^2(\Z, \D)$ the functions $a, \fa^{\pm}$ are outer, $a(0) > 0$. The reflection coefficient $\rc_{q} = b/a$ of $q$ belongs to the unit ball of $L^\infty(\T)$. It is completely determined by $b$, and, conversely, it determines the pair $a$, $b$ uniquely.
\end{Prop}
\beginpf By definition and \eqref{eq: ratio of Szego functions equals schur }, we have
\begin{gather}\label{eq: a as the product of outers}
	\fa^{\pm} = \frac{1}{2}D_{\mu}^{-1} \cdot\left(1 + \frac{1+zf^\pm}{1-zf^\pm}\right),\qquad a = \fa^+\fa^-\left(1-\fb^+\fb^-/\fa^+\fa^-\right) = \fa^+\fa^-\left(1-f^+ f^-\right). 
\end{gather}
We know that $\tfrac{1}{2}D_{\mu}^{-1}$ is outer, $1 + \frac{1+zf^\pm}{1-zf^\pm}$ and $1-f^+ f^-$ are analytic in $\D$ and have positive real part hence they are also outer, see Corollary 4.8 in \cite{Garnett}. Therefore  $\fa^{\pm}$, $a$ are outer as the products of outer functions. Next, $D_{\mu_d^{\pm}}(0)=D_{\mu^{\pm}}(0) > 0$ hence $\fa^\pm(0)$ are real and positive. We have $\beta_0 = 0$, therefore $f^{-}(0) = 0$ (recall Schur's algorithm \eqref{eq20}) and $\fb^-(0) = 0$. Thus $a(0) = \fa^+(0)\fa^-(0) > 0$. From \eqref{eq32} we have
\begin{equation}\label{eq2new}
	1 - |\rc_{q}|^2 = \frac{|a|^2 - |b|^2}{|a|^2} = \frac{1}{|a|^2} \ge 0
\end{equation}
almost everywhere on $\T$. In particular, $\rc_q$ belongs to the unit ball of $L^\infty(\T)$. We proved that $a$ is outer hence it is completely defined by $|a|$. 
Therefore, knowing the coefficient $b$, one can recover $|a| =\sqrt{ 1 + |b|^2}$ and $a$.
In particular, the numerator $b$ determines the whole fraction $\rc_q = b/a$. Conversely, if the function $\rc_q$ is given, then $|a|$ is defined by \eqref{eq2new}, hence the pair $a$, $b$ could be found from the fraction $\rc_q = b/a$. \qed

\medskip

 Next proposition shows that \eqref{eq25bis} has sense for all $q \in \ell^2(\Z, \D)$, and, moreover, the definitions of $a$, $b$ in \eqref{eq26}, \eqref{eq25bis} are equivalent. 
\begin{Prop}\label{pr1}
For every $q \in \ell^2(\Z, \D)$, the product in \eqref{eq25bis} converges in Lebesgue measure on $\T$. Moreover, the functions $a$, $b$ in \eqref{eq25bis} coincide with those in \eqref{eq26}.
\end{Prop}
\beginpf Denote by $\Phi_{\pm,n}$, $\Psi_{\pm,n}$ the monic orthogonal polynomials of $\mu^{\pm}$ and $\mu^{\pm}_d$, and let $\phi_{\pm,n}$, $\psi_{\pm,n}$ be the corresponding normalized polynomials, see \eqref{eq31}. For each $n \ge 0$, $z \in \T$, we have
\begin{gather}\notag
\prod_{k = 0}^{n}
\frac{1}{\sqrt{1-|\alpha_k|^2}}
\cdot
\prod_{k = 0}^{n}
\begin{pmatrix}
1 & \ov{\alpha_k} \bar z^k\\
\alpha_k z^k & 1
\end{pmatrix} =
\begin{pmatrix}
\frac{\psi_{+,n+1}^* + \phi_{+,n+1}^*}{2}
 & \ov{\frac{\psi^{*}_{+,n+1} - \phi^{*}_{+,n+1}}{2z}}\\
\frac{\psi^{*}_{+,n+1} - \phi^{*}_{+,n+1}}{2z} & \ov{\frac{\psi_{+,n+1}^* + \phi_{+,n+1}^*}{2}}
\end{pmatrix}.
\end{gather}
The proof is a routine verification of the identity
$$
\begin{pmatrix}
\frac{\psi_{+,n+1}^* + \phi_{+,n+1}^*}{2}
 & \ov{\frac{\psi^{*}_{+,n+1} - \phi^{*}_{+,n+1}}{2z}}\\
\frac{\psi^{*}_{+,n+1} - \phi^{*}_{+,n+1}}{2z} & \ov{\frac{\psi_{+,n+1}^* + \phi_{+,n+1}^*}{2}}
\end{pmatrix}
=
\begin{pmatrix}
\frac{\psi_{+,n}^* + \phi_{+,n}^*}{2}
 & \ov{\frac{\psi^{*}_{+,n} - \phi^{*}_{+,n}}{2z}}\\
\frac{\psi^{*}_{+,n} - \phi^{*}_{+,n}}{2z} & \ov{\frac{\psi_{+,n}^* + \phi_{+,n}^*}{2}}
\end{pmatrix}
\frac{1}{\sqrt{1-|\alpha_n|^2}}
\cdot
\begin{pmatrix}
1 & \ov{\alpha_n} \bar z^n\\
\alpha_n z^n & 1
\end{pmatrix},
$$ 
using relations \eqref{recurrence relation} and \eqref{eq43}. 
It is known that $\phi_{\pm,n}^* \to D_{\mu^{\pm}}^{-1}$, $\psi_{\pm,n}^* \to D_{\mu_{d}^{\pm}}^{-1}$ in Lebesgue measure on~$\T$, see (2.4.34) in \cite{Simonbook1}. Therefore, we have
\begin{equation}\label{eq46}
\prod_{k = 0}^{\infty}
\frac{1}{\sqrt{1-|q_k|^2}}
\cdot
\prod_{k = 0}^{\infty}
\begin{pmatrix}
	1 & \ov{q_k} \bar z^k\\
	q_k z^k & 1
\end{pmatrix} =
\prod_{k = 0}^{\infty}
\frac{1}{\sqrt{1-|\alpha_k|^2}}
\cdot
\prod_{k = 0}^{\infty}
\begin{pmatrix}
1 & \ov{\alpha_k} \bar z^k\\
\alpha_k z^k & 1
\end{pmatrix} =
\begin{pmatrix}
\fa^+ 
 & \ov{\fb^+}\\
\fb^+ & \ov{\fa^+}
\end{pmatrix},
\end{equation}
where the product converges in Lebesgue measure on $\T$.
Recall that $\beta_k = -\ov{q_{-k}}$ for $k \ge 1$, $\beta_0 = 0$. We have
\begin{align*}
\left(\prod_{k = -n}^{-1}
\frac{1}{\sqrt{1-|q_k|^2}}
\cdot
\prod_{k = -n}^{-1}
\begin{pmatrix}
1 & \ov{q_k} \bar z^k\\
q_k z^k & 1
\end{pmatrix} \right)^{-1}
&=
\prod_{k = 0}^{n}
\left(
\frac{1}{\sqrt{1-|\beta_k|^2}}
\begin{pmatrix}
1 & -\beta_k z^k\\
-\ov{\beta_k} \bar z^k & 1
\end{pmatrix}\right)^{-1}\\
&=\prod_{k = 0}^{n}
\frac{1}{\sqrt{1-|\beta_k|^2}}
\cdot
\prod_{k = 0}^{n}
\begin{pmatrix}
1 & \beta_k z^k\\
\ov{\beta_k} \bar z^k & 1
\end{pmatrix}.
\end{align*}
Note that for each $k \ge 0$ we have
$$
j_0
\begin{pmatrix}
1 & \beta_k z^k\\
\ov{\beta_k} \bar z^k & 1
\end{pmatrix}
j_0
=
\begin{pmatrix}
1 & \ov{\beta_k} \bar z^k\\
\beta_k z^k & 1
\end{pmatrix},
\qquad j_0 = 
\begin{pmatrix}
0&1\\
1&0\\
\end{pmatrix},
$$
and $\{\beta_k\}_{k \ge 0}$ coincides with the sequence of recurrence coefficients of $\mu^{-}$. 
So, we obtain
$$
j_0\left(\prod_{k = -n}^{-1}
\frac{1}{\sqrt{1-|q_k|^2}}
\cdot
\prod_{k = -n}^{-1}
\begin{pmatrix}
1 & \ov{q_k} \bar z^k\\
q_k z^k & 1
\end{pmatrix} \right)^{-1} j_0 \to 
\begin{pmatrix}
\fa^- & \ov{\fb^-}\\
\fb^- & \ov{\fa^-}
\end{pmatrix},
$$
where the convergence is in Lebesgue measure on $\T$. Taking the inverses (note that $|\fa^-|^2 - |\fb^-|^2 = 1$ from the consideration of determinants), we obtain
$$
\prod_{k = -\infty}^{-1}
\frac{1}{\sqrt{1-|q_k|^2}}
\cdot
\prod_{k = -\infty}^{-1}
\begin{pmatrix}
1 & \ov{q_k} \bar z^k\\
q_k z^k & 1
\end{pmatrix} = j_0 \begin{pmatrix}
\ov{\fa^-} & -\ov{\fb^-}\\
-\fb^- & \fa^-
\end{pmatrix}j_0 = \begin{pmatrix}
\fa^- & -\fb^-\\
-\ov{\fb^-} & \ov{\fa^-}
\end{pmatrix}.
$$
Eventually, we get
\begin{equation}\label{eq27}
\prod_{k = -\infty}^{\infty}
\frac{1}{\sqrt{1-|q_k|^2}}
\cdot
\prod_{k = -\infty}^{\infty}
\begin{pmatrix}
1 & \ov{q_k} \bar z^k\\
q_k z^k & 1
\end{pmatrix}
= 
\begin{pmatrix}
\fa^- & -\fb^-\\
-\ov{\fb^-} & \ov{\fa^-}
\end{pmatrix}
\begin{pmatrix}
\fa^+ 
 & \ov{\fb^+}\\
\fb^+ & \ov{\fa^+}
\end{pmatrix}
=
\begin{pmatrix}
a & b\\
\ov{b} & \bar a
\end{pmatrix},
\end{equation}
with $a = \fa^+\fa^{-} - \fb^+\fb^-$, $b = \fa^-\ov{\fb^+} - \fb^-\ov{\fa^+}$, as claimed. \qed

\medskip

We can now prove Lemma \ref{l34} from Section \ref{s3}.

\medskip

\noindent {\bf Proof of Lemma \ref{l34}.}\label{pfl34} Propositions \ref{pr1}, \ref{pr2} imply that the definitions of $a$, $b$ in \eqref{eq25bisbis} and \eqref{eq26} are equivalent. Note that for $q \in \ell^2(\Z, \D)$ with $\supp q \subset \Z_+$ we have $\fa^- = 1, \fb^-=0$ hence
$$
\fc_{q} = \frac{\bar b}{a} = \frac{\fb_+}{\fa_+} = f^+. 
$$
In particular, the recurrence coefficients of $\fc_{q}$ coincide with those of $f^+$, $\mu^+$, i.e., with the sequence $\{q(k)\}_{k \in \Z_+}$. \qed

\medskip 

\begin{Prop}\label{pr9}
We have $\rc_{q(\cdot - n)} = z^{-n}\rc_{q}$ for every compactly supported $q \in \ell^2(\Z, \D)$ and $n\in \Z$.
\end{Prop}
\beginpf We have
\begin{align*}
\prod_{k \in \Z}\begin{pmatrix}
1 & \ov{q(k-n)}z^{-k}\\
q(k-n)z^{k} & 1
\end{pmatrix}
&=
\prod_{k \in \Z}
\begin{pmatrix}
1 & 0\\
0 & z^{n}
\end{pmatrix}
\begin{pmatrix}
1 & \ov{q(k-n)}z^{-(k-n)}\\
q(k-n)z^{k-n} & 1
\end{pmatrix}
\begin{pmatrix}
1 & 0\\
0 & z^{-n}
\end{pmatrix}\\
&=
\begin{pmatrix}
1 & 0\\
0 & z^{n}
\end{pmatrix}
\left[\prod_{k \in \Z}
\begin{pmatrix}
1 & \ov{q(k-n)}z^{-(k-n)}\\
q(k-n)z^{k-n} & 1
\end{pmatrix}\right]
\begin{pmatrix}
1 & 0\\
0 & z^{-n}
\end{pmatrix}\\
&=\begin{pmatrix}
1 & 0\\
0 & z^{n}
\end{pmatrix}
\begin{pmatrix}
a & b\\
\bar b & \bar a
\end{pmatrix}
\begin{pmatrix}
1 & 0\\
0 & z^{-n}
\end{pmatrix}
=\begin{pmatrix}
a & bz^{-n}\\
\bar b z^{n} & \bar a
\end{pmatrix}.
\end{align*}
Hence $\rc_{q(\cdot - n)} = b z^{-n}/a = z^{-n}\rc_{q}$ by \eqref{eq32} and  Proposition \ref{pr1}. \qed

\medskip

\begin{Prop}\label{pr11}
There are $q_1 \neq q_2$ in $\ell^2(\Z, \D)$ such that $\rc_{q_1} = \rc_{q_2}$.
\end{Prop}
\beginpf Following \cite{TT}, let us consider an imaginary-valued function $\fb$ on $\T$ of Smirnov class in the unit disk. One can take, say, $\fb = \frac{1+z}{1-z}$. Let $\fa$ be the outer function in $\D$ such that $|\fa|^2 - |\fb|^2 = 1$ almost everywhere on $\T$. The function $f = \fb/\fa$ is a Schur function of Szeg\H{o} class. Indeed, $\log(1-|f|^2)=\log|a|^{-2}$ belongs to $L^1(\T)$. Therefore, we can define the sequences  $q_1,\tilde q_2\in \ell^2(\Z, \D)$ by
\begin{gather*}
\begin{matrix}
n: \phantom{=}&\ldots & -3 & -2 & -1 & 0 & 1 & 2 &\ldots&
\\
q_1 =(&\ldots & 0, & 0, & 0, & f_{0}(0), & f_1(0), & f_2(0), &  \ldots&)
\\
\tilde q_2 = (&\ldots & -\ol{f_2(0)}, & -\ol{f_1(0)}, & -\ol{f_0(0)}, & 0, & 0, & 0, &\ldots&).
\end{matrix}
\end{gather*}
For these sequences, we  have
\begin{gather*}
\fa^{+}_{q_1} = \fa, \;\; \fb^{+}_{q_1} = \fb,\qquad \fa^-_{q_1} = 1,\;\; \fb_{q_1}^{-} = 0, \qquad \fa^{+}_{\tilde q_2} = 1,\;\; \fb^{+}_{\tilde q_2} = 0. 
\end{gather*}
Furthermore, from the proof of Proposition \ref{pr9} we obtain $\fb_{\tilde q_2}^{-} = z\fb$ and $\fa^-_{\tilde q_2} = \fa$. Therefore,
\begin{align*}
&a_{q_1} = \fa \cdot 1 - \fb \cdot 0 = \fa,  &b_{q_1} = 1 \cdot \ov{\fb} - 0 \cdot \ov{\fa} = \ov{\fb},\\
&a_{\tilde q_2} = 1 \cdot \fa - 0 \cdot z\fb = \fa,  &b_{\tilde q_2} = \fa \cdot \ov{0} - z\fb \cdot \ov{1} = -z\fb.
\end{align*}
Then $\rc_{q_1} = \frac{\ov{\fb}}{\fa}$, $\rc_{\tilde q_2} = \frac{-z\fb}{\fa}$, and, since $\fb = -\bar\fb$, we have $\rc_{q_1} = z\rc_{\tilde q_2}$ almost everywhere on $\T$. Note that
$z\rc_{\tilde q_2} = \rc_{\tilde q_2(\cdot +1)}$ by Proposition \ref{pr9}. Now set $q_2 = \tilde q_2(\cdot +1)$ and observe that $\rc_{q_1} = \rc_{q_2}$, while $q_1$, $q_2$ are supported on disjoint subsets of $\Z$, so $q_1 \neq q_2$. \qed

\medskip

\begin{Prop}\label{pr3}
For every $q \in \ell^2(\Z, \D)$, we have 
\begin{align}\label{eq7}
\int_{\T}\log(1 - |\rc_{q}|^2)\,dm = 
-\log|a(0)|^2 =
\log\prod_{n\in\Z}(1-|q(n)|^2).
\end{align}
\end{Prop}
\beginpf Take a sequence $q$ in $\ell^2(\Z, \D)$ and define $\{\alpha_n\}$, $\{\beta_n\}$, $\mu^\pm$, $\fa^\pm$, $\fb^\pm$, $f^\pm$, $a$ and $b$ as in the beginning of Section \ref{s52}. From \eqref{eq2new} and the mean value theorem, we get
\begin{gather*}
	\int_{\T}\log(1 - |\rc_{q}|^2)\,dm 
	= -\int_{\T}\log|a|^2\,dm = -\log|a(0)|^2.
\end{gather*}
In the proof of Proposition \ref{pr2} we established $a(0) = \fa^+(0)\fa^-(0) = D_{\mu^+}^{-1}(0)D_{\mu^-}^{-1}(0)$. Let $w^\pm$ be the densities of the a.\,c.\ parts $\mu^\pm$ with respect to the Lebesgue measure on $\T$, then from formula \eqref{D in z formula} and Szeg\H{o} theorem \eqref{eq3} it follows that
\begin{align*}
	-\log a(0)^2 &= \int_{\T}\log w^+(\xi)\,dm(\xi) + \int_{\T}\log w^-(\xi)\,dm(\xi)  
	\\
	&= \log\prod_{n \ge 0}(1 - |\alpha_n|^2) + \log\prod_{n \ge 0}(1 - |\beta_n|^2) = \log\prod_{n \in \Z}(1 - |q_n|^2),
\end{align*}
as claimed. \qed
\medskip

\begin{Prop}\label{pr5}
For every $q \in \ell^2(\Z, \D)$, the functions $\fa^{\pm}/{a}$, $\fb^{\pm}/{a}$ belong to the unit ball of the Hardy class $H^2(\D)$.
\end{Prop}
\beginpf Since $a$, $\fa^{\pm}$ are outer in $\D$ and $\fb^{\pm}$ are in the Smirnov class (see \eqref{eq49}), we need to show only that $\fa^{\pm}/{a} $ belong to the unit ball of $L^2(\T)$. Denote, as before, $f^{\pm} = \fb^{\pm}/\fa^{\pm}$, and recall that $f^{\pm}$ are Schur functions. The function 
$$
h = \frac{1-|f^{-}f^{+}|^2}{|1 - f^{-}f^{+}|^2}
=
\Re\left(\frac{1 + f^{-}f^{+}}{1 - f^{-}f^+}\right) 
$$
is positive and harmonic in $\D$, therefore, it coincides with the Poisson integral of a finite positive Borel measure on $\T$. Moreover, $h$ is equal to the density of the absolutely continuous part of that measure almost everywhere on $\T$. Hence, $h \in L^1(\T)$ (we borrowed this trick from \cite{TT}) and 
\begin{gather*}
	\|h\|_{L^1(\T)} = \int_{\T}\frac{1-|f^{-}f^{+}|^2}{|1 - f^{-}f^{+}|^2}\,dm = \Re\left(\frac{1 + f^{-}(0)f^{+}(0)}{1 - f^{-}(0)f^+(0)}\right) = 1,
\end{gather*}
because $f^{-}(0) = 0$. On the other hand, by \eqref{eq: a as the product of outers} we have
\begin{gather*}
	\frac{1}{|\fa^{\pm}|^2} = 1 - \frac{|\fb^{\pm}|^2}{|\fa^{\pm}|^2} = 1 - |f^\pm|^2, \qquad \frac{1}{|a|^2} = \frac{(1-|f^+|^2)(1-|f^-|^2)}{|1- f^+f^-|^2}.
\end{gather*}
almost everywhere on $\T$. It gives us 
$$
\left|\frac{\fa^{\pm}}{a}\right|^2 = \frac{(1-|f^+|^2)(1-|f^-|^2)}{(1-|f^\pm|^2)|1- f^+f^-|^2} = \frac{1-|f^{\mp}|^2}{|1- f^+f^-|^2}
\le 
\frac{1-|f^{-}f^{+}|^2}{|1 - f^{-}f^{+}|^2} = h.
$$
Therefore $\|\fa^{\pm}/a\|_{L^2(\T)}^{2}\le \|h\|_{L^1(\T)} = 1$, as claimed. \qed

\medskip

The authors are grateful to S.~Denisov for the argument based on \eqref{eq11} in the proof of proposition below.
\begin{Prop}\label{pr4}
Suppose that $q_1, q_2 \in \ell^2(\Z, \D)$ are such that $\rc_{q_1} = \rc_{q_2}$. If 
$\|\rc_{q_{1,2}}\|_{L^\infty(\T)} < 1$, 
then $q_1 = q_2$.
\end{Prop}
\beginpf Let $a_{1,2}$, $b_{1,2}$ be the coefficients in \eqref{eq26} corresponding to $\rc_{q_{1}}$, $\rc_{q_{2}}$, respectively. 
By Proposition \ref{pr2}, we have $a_1 = a_2$, $b_1 = b_2$, so we denote $a=a_{1,2}$, $b = b_{1,2}$. Then \eqref{eq26} gives four identities
$$
a = \fa_{k}^+\fa_{k}^{-} - \fb_{k}^+\fb_{k}^-, \qquad 
b = \fa_{k}^{-}\ov{\fb_{k}^+} - \fb_{k}^{-}\ov{\fa_{k}^+}, \qquad k =1,2,
$$
for the functions $\fa_{1,2}^{\pm}$, $\fb_{1,2}^{\pm}$ corresponding to $q_1$, $q_2$.
A simple algebra yields
\begin{equation}\label{RH}
\begin{pmatrix}
\fa_{k}^+ & \fb_{k}^{+} \\
\fb_{k}^{-} & 
\fa_{k}^{-} 
\end{pmatrix} 
\overline{
\begin{pmatrix}
\fa_{k}^+ & -\fb_{k}^{-} \\
-\fb_{k}^{+} & \fa_{k}^{-}
\end{pmatrix}
}
=
\begin{pmatrix}
|\fa^+_{k}|^2 - |\fb^+_{k}|^2 & -\fa^+_{k}\ov{\fb^{-}_{k}} + \fb_{k}^{+}\ov{\fa_{k}^{-}} \\
\fb^{-}_{k}\ov{\fa^+_{k}} - \fa_{k}^{-}\ov{\fb_{k}^{+}} & |\fa^-_{k}|^2 - |\fb^-_{k}|^2
\end{pmatrix}
=
\begin{pmatrix}
1 & \ov{b} \\
-b & 1
\end{pmatrix}, 
\end{equation}
almost everywhere on $\T$ for $k = 1,2$. In particular, we have
$$
\begin{pmatrix}
\fa_{1}^+ & \fb_{1}^{+} \\
\fb_{1}^{-} & 
\fa_{1}^{-} 
\end{pmatrix} 
\overline{
\begin{pmatrix}
\fa_{1}^+ & -\fb_{1}^{-} \\
-\fb_{1}^{+} & \fa_{1}^{-}
\end{pmatrix}
}
=
\begin{pmatrix}
\fa_{2}^+ & \fb_{2}^{+} \\
\fb_{2}^{-} & 
\fa_{2}^{-} 
\end{pmatrix} 
\overline{
\begin{pmatrix}
\fa_{2}^+ & -\fb_{2}^{-} \\
-\fb_{2}^{+} & \fa_{2}^{-}
\end{pmatrix}
}.
$$
Inverting matrices in the last equation, we obtain
\begin{align}
	\notag
\begin{pmatrix}
	\fa_{2}^+ & \fb_{2}^{+} \\
	\fb_{2}^{-} & 
	\fa_{2}^{-} 
\end{pmatrix} ^{-1}
\begin{pmatrix}
	\fa_{1}^+ & \fb_{1}^{+} \\
	\fb_{1}^{-} & 
	\fa_{1}^{-} 
\end{pmatrix} 
=
\overline{
\begin{pmatrix}
	\fa_{2}^+ & -\fb_{2}^{-} \\
	-\fb_{2}^{+} & \fa_{2}^{-}
\end{pmatrix}
\begin{pmatrix}
	\fa_{1}^+ & -\fb_{1}^{-} \\
	-\fb_{1}^{+} & \fa_{1}^{-}
\end{pmatrix}^{-1}
},
\\
\label{eq11}
I:=\frac{1}{a}
\begin{pmatrix}
\fa_{2}^{-} & 
-\fb_{2}^{+} \\
-\fb_{2}^{-} & 
\fa_{2}^+ 
\end{pmatrix}
\begin{pmatrix}
\fa_{1}^+ & 
\fb_{1}^{+} \\
\fb_{1}^{-} & 
\fa_{1}^{-} 
\end{pmatrix} 
= 
\overline{
\frac{1}{a}
\begin{pmatrix}
\fa_{2}^+ & -\fb_{2}^{-} \\
-\fb_{2}^{+} & \fa_{2}^{-}
\end{pmatrix}
}
\overline{
\begin{pmatrix}
\fa_{1}^{-} & 
\fb_{1}^{-} \\
\fb_{1}^{+} & 
\fa_{1}^+
\end{pmatrix}
}.
\end{align} 
Equating the $(1,1)$ matrix elements in this identity, we get
\begin{gather*}
\frac{\fa^+_1\fa^-_2 - \fb_1^{-}\fb_2^{+}}{a} = \overline{\left(\frac{\fa^+_2\fa^-_1 - \fb_2^{-}\fb_1^{+}}{a}\right)}.
\end{gather*}
Formula \eqref{eq2new} and our assumption $\|\rc_{q_{1,2}}\|_{L^\infty(\T)} < 1$ imply that $a \in H^\infty(\D)$. We now see from Proposition \ref{pr5} that the functions $F_1 = \frac{\fa^+_1\fa^-_2 - \fb_1^{-}\fb_2^{+}}{a}$, $F_2 = \frac{\fa^+_2\fa^-_1 - \fb_2^{-}\fb_1^{+}}{a}$ belong to the Hardy space $H^{1}(\D)$. Therefore $F_1$ and $F_2$ are constant functions and 
$$
\ov{F_2} = F_1 = F_1(0) = \frac{\fa^+_1(0)\fa^-_2(0) - \fb_1^{-}(0)\fb_2^{+}(0)}{a(0)} =  \frac{a(0)}{a(0)} =1.
$$ 
In other words, the $(1,1)$ coefficient of the matrix $I$ in \eqref{eq11} is $1$. Note that it coincides with the $(2,2)$ coefficient of $I$. Similarly, we use $\fb_1^{\pm}(0)=\fb_2^{\pm}(0) = 0$ and prove that $(1,2)$, $(2,1)$ coefficients of $I$ are $0$ thus getting
$$
\frac{1}{a}
\begin{pmatrix}
\fa_{2}^{-} & 
-\fb_{2}^{+} \\
-\fb_{2}^{-} & 
\fa_{2}^+ 
\end{pmatrix}
\begin{pmatrix}
\fa_{1}^+ & 
\fb_{1}^{+} \\
\fb_{1}^{-} & 
\fa_{1}^{-} 
\end{pmatrix} = \begin{pmatrix}1 & 0\\ 0 & 1\end{pmatrix},
$$
which is equivalent to
$$
\begin{pmatrix}
\fa_{1}^+ & 
\fb_{1}^{+} \\
\fb_{1}^{-} & 
\fa_{1}^{-} 
\end{pmatrix} = 
\begin{pmatrix}
\fa_{2}^+ & \fb_{2}^{+}\\
\fb_{2}^{-} & \fa_{2}^-
\end{pmatrix}.
$$
It follows that $f^{\pm}_{1} = f^{\pm}_{2}$, which, in turn, is equivalent to $q_1 = q_2$ on $\Z$, because the recurrence coefficients of $f_{1,2}^{\pm}$ determine completely $q_{1,2}$ on $\Z_{\pm}$, see the beginning of Section \ref{s52}. \qed

\subsection{Convergence in the space $X$}
We first prove a version of Sylvester--Winebrenner theorem \cite{Sylvester} for Schur functions. Let us recall its statement.
\begin{Prop}[Sylvester--Winebrenner theorem]\label{pr6} The mapping 
$f \mapsto \{f_n(0)\}_{n \ge 0}$ that takes a Schur function into the sequence of its recurrence coefficients is a homeomorphism from the metric space $X_+ = \{f \in \Sch_{*}(\D): \eta(F) >0 \}$ 
with the metric $\rho_s(f,g)^2 = -\int_{\T}\log\bigl(1 - \bigl|\frac{f-g}{1-\bar f g}\bigr|^2\bigr)\,dm$ onto the metric space $\ell^2(\Z_+, \D)$ of sequences $q: \Z_+ \to \D$ with the metric 
$\|q - \tilde q\|_{\ell^2}^{2} = \sum_{n \in \Z_+}|q(n) - \tilde q(n)|^2$. 
\end{Prop}
\beginpf Assume that $f_n, f \in X_+$ are such that $\rho_s(f_n,f) \to 0$. Let $q_n$, $q$ be the sequences of recurrence coefficients of $f_n$, $f$, respectively. By Szeg\H{o} theorem, we have $q_n, q \in \ell^2(\Z_+, \D)$, and, moreover, 
$$
-\log \prod_{k\ge 0}(1-|q_n(k)|^2) = \rho_s(f_n,0) \to \rho_s(f,0) = -\log \prod_{k\ge 0}(1-|q(k)|^2).
$$
The convergence $q_n \to q$ in $\ell^2(\Z_+, \D)$ will follow if we check that $q_n(k) \to q(k)$ for each $k \in \Z_+$ (indeed, we then have $\sum_{k \ge N}|q_n(k)|^2 \to 0$ as $N \to +\infty$ uniformly in $n \in \Z_+$). To this end, note that assumption $\rho_s(f_n,f) \to 0$ implies that the sequence $\{f_n\}$ converges to $f$ in Lebesgue measure on $\T$, and, since $|f_n| \le 1$, $|f| \le 1$ on $\T$, the functions $f_n$ converge to $f$ uniformly on compacts in $\D$. Now the fact that $q_n(k) = (f_n)_k(0)$ tends to $(f)_k(0) = q(k)$ as $n \to +\infty$ for every $k \in \Z_+$ follows from Schur's algorithm \eqref{eq: NFLT schur algorithm}. We see that the mapping $f \mapsto q$ is continuous from $X_+$ to $\ell^2(\Z_+, \D)$.

\smallskip 

\noindent Turning to the inverse mapping, we introduce the quantities (see \cite{Sylvester})
\begin{gather}\label{eq: E defininiton}
	E(f, g) = -\int_{\T}\log(1 - \bar f g)\,dm, \qquad E(f) = E(f,f).
\end{gather}
We have $1 - \bigl|\frac{f-g}{1-\bar f g}\bigr|^2 = \frac{(1 - |f|^2)\cdot (1 -|g|^2)}{|1-\bar f g|^2} $, hence
\begin{gather}\label{eq47}
\rho_s(f,g)^2 = E(f) + E(g) - 2\Re E(f,g).
\end{gather}
Suppose that $q_n, q$ are sequences in $\ell^2(\Z_+,\D)$ such that $q_n \to q$ in $\ell^2(\Z_+,\D)$. Denote by $f_n$, $f$ the Schur functions corresponding to these sequences. We have $f_n, f \in X_+$ by Szeg\H{o} theorem, see \eqref{eq3}. Let us prove that $\rho_{s}(f_n,f) \to 0$ as $n\to +\infty$. Since
$E(f_n) \to E(f)$ by Szeg\H{o} theorem, relation \eqref{eq47} shows that we only need to check that $E(f_n, f) \to E(f, f)$. We have 
$$
E(f_n, f) = -\int_{\T}\sum_{k \ge 0}\frac{(\bar f_n f)^k}{k}\,dm = -\int_{\T}\sum_{k=1}^{N}\frac{(\bar f_n f)^k}{k}\,dm -
\int_{\T}\sum_{k=N+1}^{\infty}\frac{(\bar f_n f)^k}{k}\,dm,
$$
and 
$$
\left|\int_{\T}\sum_{k=N+1}^{\infty}\frac{(\bar f_n f)^k}{k}\,dm\right| \le \int_{\T}\sum_{k=N+1}^{\infty}\frac{|\bar f_n f|^k}{k}\,dm \le \int_{\T}\sum_{k=N+1}^{\infty}\frac{|f|^k}{k}\,dm, 
$$
which tends to zero as $N \to +\infty$ by Lebesgue dominated convergence theorem (the majorant is $\log\frac{1}{1-|f|} \in L^1(\T)$). Next, let us show that for each $k\in \Z_+$ we have
\begin{gather}\label{from f_n^k to f^k}
	\int_{\T}(\bar f_n f)^k\,dm \to \int_{\T}|f|^{2k}\,dm, \qquad n \to +\infty.
\end{gather}
Indeed, the first $m$ Taylor coefficients of $f$ are polynomials in $q(0),\ov{q(0)}, \ldots, q(m-1), \ov{q(m-1)}$ and similarly for $f_n$, see Lemma \ref{lemma: schur functions witl the same an} of Section 1.3 in \cite{Simonbook1}. Hence Taylor coefficients of $f_n^k$ tend to those of $f^k$ as $n\to\infty$. Rewrite quantity in \eqref{from f_n^k to f^k} as
\begin{gather*}
	\int_{\T}(\bar f_n f)^k\,dm = \sum_{m = 0}^{\infty}\ol{c_m(f_n^k)}c_m(f^k) = \sum_{m = 0}^{M}\ol{c_m(f_n^k)}c_m(f^k) + \sum_{m = M + 1}^{\infty}\ol{c_m(f_n^k)}c_m(f^k).
\end{gather*}
The second sum can be estimated using the Cauchy inequality by 
\begin{gather*}
	\|f_n^k\|_{H^2(\D)}^2\cdot\left(\sum_{m = M + 1}^{\infty}|c_m(f^k)|^2\right) \le \sum_{m = M + 1}^{\infty}|c_m(f^k)|^2,
\end{gather*}
because $f_n^k \in \Sch(\D)$ and consequently $\|f_n^k\|_{H^2(\D)}^2\le 1$. Hence it tends to $0$ as $M\to\infty$. The first sum  tends to $\sum_{m = 0}^{M}|c_m(f^k)|^2$ as $n\to\infty$ and \eqref{from f_n^k to f^k} follows. Relation \eqref{from f_n^k to f^k} shows that $E(f_n, f) \to 0$, $\rho_s(f_n,f) \to 0$, and thus the mapping $q \mapsto f$ is continuous from $\ell^2(\Z_+, \D)$ to the metric space $X_+$. \qed

\medskip

The following lemma is elementary. It is known as Scheff\'e's lemma, see Section 5.10 in \cite{williams1991probability}.

\begin{Lem}\label{l5}
Let measurable functions $g$, $g_n$ on $\T$ be such that $g_j\to g$ in Lebesgue measure on $\T$ and $\|g_j\|_{L^1(\T)}\to \|g\|_{L^1(\T)} $ as $j\to\infty$. Then $\|g-g_j\|_{L^1(\T)}\to 0$.
\end{Lem}
\beginpf If $\|g\|_{L^1(\T)} = 0$, then the lemma is trivial, otherwise we can reduce the statement of the lemma to the case $\|g_j\|_{L^1(\T)}= \|g\|_{L^1(\T)}  = 1$ by changing $g$ and $g_n$ to $g/\|g\|_{L^1(\T)}$ and $g_n/\|g_n\|_{L^1(\T)}$ respectively. Consider any subsequence $g_{n_k}$ of the sequence $g_n$. Let $g_{n_{k_j}}$ be its subsequence converging Lebesgue almost everywhere on $\T$. The limit of $g_{n_{k_j}}$ coincides with $g$ Lebesgue almost everywhere on~$\T$. To simplify notation, we denote the new sequence $g_{n_{k_j}}$ by $\tilde g_j$. Let $\eps>0$. By Egorov's theorem and integrability of $g$, there is $K_\eps\subset \T$ such that $m(K_\eps)<\eps, \|g\|_{L^1(K_\eps)}<\eps$ and $\tilde  g_j\to g$ uniformly on $\T\setminus K_\eps$. In particular, we have
\[
\int_{\T\setminus K_\eps} |\tilde g_j|\, dm \to \int_{\T\setminus K_\eps}|g|\,dm \ge 1 -2\eps,
\qquad
\limsup_{j\to\infty}\int_{K_\eps} |\tilde g_j|\, dm \le 2\eps.
\]
Now, we only need to write
\[
\limsup_{j\to\infty}\|g-\tilde  g_j\|_{L^1(\T)}\le \limsup_{j\to\infty}\|g-\tilde  g_j\|_{L^1(\T\setminus K_\eps)}+\limsup_{j\to\infty}\|\tilde  g_j\|_{L^1(K_\eps)}+\|g\|_{L^1(K_\eps)} \le 3\eps.
\]
Since $\eps>0$ is arbitrary, we see that $\tilde  g_j \to g$ in $L^1(\T)$. In other words, we have shown that any subsequence of $g_n$ contains a subsequence converging to $g$ in $L^1(\T)$. Then $g_n \to g$ in $L^1(\T)$ and the lemma follows. \qed

\medskip

Recall that the space $X$ and the metric $\rho_s$ on $X$ are defined in \eqref{eq: X space def} and \eqref{eq: sylvester metric definition}. For $r \in X$, define the function $E(r)$ by \eqref{eq: E defininiton}.

\begin{Prop}\label{prop: convergence in measure to convergence in X}
Let $r_n,r \in X$. The following assertions are equivalent:
\begin{itemize}
	\item[$(a)$] $r_n$ converges to $r$ in $X$;
	\item[$(b)$]  $r_n$ converges to $r$ in Lebesgue measure on $\T$ and $\displaystyle\lim_{n\to+\infty}E(r_n)= E(r)$.
\end{itemize}
\end{Prop}
\beginpf Assume that $r_n \to r$ in $X$ as $n \to +\infty$. The convergence in measure follows immediately. For all $n\ge 0$, we have $|1-\ov{r_n}r| \ge 1 - |r|$ and $\log\frac{1}{1-|r|} \in L^1(\T)$. Hence by the dominated convergence theorem we have
\begin{gather}\label{eq: E(r_n, r) to E(r)}
	E(r_n, r) = -\int_{\T}\log (1-\ov{r_n}r)\,dm 
	\to 
	-\int_{\T}\log (1-|r|^2)\,dm = E(r).
\end{gather}
Thus, from \eqref{eq47} we see that 
$$
 0 = \lim_{n \to +\infty} \rho_{s}(r_n, r)^2 = \lim_{n \to +\infty}(E(r) + E(r_n) - 2\Re E(r, r_n)) = \lim_{n \to +\infty}(E(r_n) - E(r)),
$$
which gives us the required assertion. On the other hand, if we assume $(b)$, then \eqref{eq: E(r_n, r) to E(r)} will follow by the same argument and similarly by \eqref{eq47} we will get
\begin{gather}
	\lim_{n \to +\infty} \rho_{s}(r_n, r)^2 = \lim_{n \to +\infty}(E(r) + E(r_n) - 2\Re E(r, r_n)) = 0,
\end{gather}
which is the convergence in $X$.

\qed

\begin{Prop}\label{pr7}
If $q_n \to q$ in $\ell^2(\Z, \D)$, then $\rc_{q_n} \to \rc_{q}$ in $X$.
\end{Prop}
\beginpf We want to apply the criteria from Proposition \ref{prop: convergence in measure to convergence in X}. Convergence in $\ell^2(\Z, \D)$ implies the convergence 
\begin{gather*}
	\prod_{k\in\Z}(1-|q_n(k)|^2)\to \prod_{k\in\Z}(1-|q(k)|^2),\quad n\to\infty, 
\end{gather*}
which yields $E(\rc_{q_n})\to E(\rc_q)$ by Proposition \ref{pr3}. Thus, it is suffices to show only that $\rc_{q_n} \to \rc_{q}$ in Lebesgue measure on $\T$. Recall that for every $q \in \ell^2(\Z, \D)$, we have $f^\pm = \fa^\pm / \fb^\pm$, where $f^\pm$ are the Schur functions generated by $q$, hence
\begin{gather*}
\rc_{q} 
= \frac{\fa^-\ov{\fb^+} - \fb^-\ov{\fa^+}}{\fa^+ \fa^- - \fb^+\fb^-} = \frac{\ov{\fa^{+}}}{\fa^{+}}\frac{\ov{f^+} - f^{-}}{1 - f^+f^-}\notag = \exp(-2i\Hh(\log|\fa^+|))\frac{\ov{f^+} - f^{-}}{1 - f^+f^-}.
\end{gather*}
Here $\Hh$ denotes the Hilbert transform and we used the fact that $\fa^+$ is an outer function. Furthermore we have $1/|\fa^+|^2 = 1 - |\fb^+|^2/|\fa^+|^2 = 1 - |f^+|^2$ hence
\begin{gather}\label{eq11new}
	\rc_{q} = \exp(i\Hh(\log(1-|f^+|^2)))\frac{\ov{f^+} - f^{-}}{1 - f^+f^-}.
\end{gather}
Similar formulae with Schur functions $f_{n}^{\pm}$ in place of $f^{\pm}$ hold for $q_n$. Proposition \ref{pr6} implies the convergence $f_n^\pm \to f^{\pm}$ in Lebesgue measure on $\T$. Moreover, by the Szeg\H{o} theorem, $\|1 - |f_n^+|^2\|_{L^1(\T)}\to\|1 - |f^+|^2\|_{L^1(\T)}$ hence Lemma~\ref{l5} can be applied to functions 
\begin{gather*}
	g_n = \log (1 - |f_{n}^{+}|^2), \qquad 
	g =  \log (1 - |f^{+}|^2).
\end{gather*}
It gives the convergence of  $\log (1 - |f_{n}^{+}|^2)$ to $\log (1 - |f^{+}|^2)$ in $L^1(\T)$. 
Weak continuity of the Hilbert transform $\Hh$ (see Section III.2 in \cite{Garnett}) then implies that $\exp(i\Hh(\log(1-|f_n^+|^2)))$ converges in Lebesgue measure to $\exp(i\Hh(\log(1-|f^+|^2)))$. From here and \eqref{eq11new} we see that functions $\rc_{q_n}$ converge to $\rc_{q}$ in Lebesgue measure on $\T$.   \qed

\medskip

The following proposition is not used in the proof of Theorem \ref{t5}, but it explains how instabilities may arise in Schur's algorithm.

\begin{Prop}\label{pr12}
There is $\eta > 0$ such that 
the mapping $f \mapsto \{f_n(0)\}_{n \ge 0}$ taking a Schur function $f$ into the sequence of its recurrence coefficients
is not uniformly continuous with respect to the metrics in $X_+$, $\ell^2(\Z,\D)$ on the subset of functions $f \in X_+$ satisfying $\eta(f) > \eta$.
\end{Prop}
\beginpf Take any $q \neq \tilde q$ in $\ell^2(\Z, \D)$ such that $\rc_{q} = \rc_{\tilde q}$, see Proposition \ref{pr11}. Fix $\eps > 0$ and use Proposition \ref{pr7} to find a number $N(\eps)$ such that $\rho_s(\rc_{q_N}, \rc_{q}) \le \eps$, $\rho_s(\rc_{\tilde q_N}, \rc_{\tilde q}) \le \eps$ for every $N \ge N(\eps)$, where $q_{N}(k) = q(k)$, $\tilde q_{N}(k) = \tilde q(k)$ for $k \le N-1$, and $q_{N}(k) = \tilde q_{N}(k) = 0$ for $k \ge N$. Next, shift these sequences to make them supported on $(-\infty, -1]$: define
$q_{N,s}(k) = q_N(k + N)$, $\tilde q_{N,s}(k) = \tilde q_N(k + N)$ for $k \in \Z$. Let also $q_{s} = q(\cdot + N)$, $\tilde q_{s} = \tilde q(\cdot + N)$. 
We have
$$
\rho_s(\rc_{q_{N,s}}, \rc_{\tilde q_{N,s}}) \le 
\rho_s(\rc_{q_{N,s}}, \rc_{q_s}) +
\rho_s(\rc_{q_s}, \rc_{\tilde q_s}) +
\rho_s(\rc_{\tilde q_s}, \rc_{\tilde q_{N,s}}) \le 2 \eps,
$$
because 
\begin{align*}
\rho_s(\rc_{q_{N,s}}, \rc_{q_s}) &=  \rho_s(\rc_{q_{N}}, \rc_{q}) \le \eps,\\ \rho_s(\rc_{q_s}, \rc_{\tilde q_s}) &= \rho_s(\rc_{q}, \rc_{\tilde q}) = 0,\\  \rho_s(\rc_{\tilde q_{N,s}}, \rc_{\tilde q_s}) &= \rho_s(\rc_{\tilde q_{N}}, \rc_{\tilde q}) \le \eps,
\end{align*}
by Proposition \ref{pr9} (it was proved for compactly supported $q$, but continuity in Proposition \ref{pr11} extends it to whole space $\ell^2(\Z, \D)$). 
On the other hand, $-\rc_{q_{N,s}}$, $-\rc_{\tilde q_{N,s}}$ coincide on $\T$ with Schur functions with the recurrence coefficients $\beta_N(n) = -\ov{q_{N,s}(-n)}$, $\tilde \beta_N(n) = -\ov{\tilde q_{N,s}(-n)}$, $n \ge 0$, respectively, see \eqref{eq48}, \eqref{eq26}, \eqref{eq32}. Since the sequences $\{\beta_N(n)\}_{n \ge 0}$, $\{\tilde \beta_N(n)\}_{n \ge 0}$ are uniformly separated in $\ell^2(\Z_+, \D)$ for large $N$, and $\rho_s(\rc_{q_{N,s}}, \rc_{\tilde q_{N,s}}) \le 2\eps$ for all $N \ge N(\eps)$, the mapping in the statement of proposition cannot be uniformly continuous.

\medskip

\begin{Prop}\label{pr8}
Let $q_n \in \ell^2(\Z, \D)$ be such that $\rc_{q_n} \to \rc$ in $X$ for some $\rc \in X$. Then there is a subsequence $q_{n_j}$ such that $q_{n_j} \to q$ in $\ell^2(\Z, \D)$ and 
$\rc = \rc_{q}$. 
\end{Prop}
\beginpf Since $\rc_{q_n} \to \rc$ in $X$, we know that $\rc_{q_n} \to \rc$ in Lebesgue measure on $\T$. Moreover, $E(\rc_{q_n}) \to E(\rc)$ as $n \to +\infty$ by Proposition \ref{prop: convergence in measure to convergence in X}. Hence Lemma \ref{l5} is applicable and we see that $\log(1-|\rc_{q_n}|^2)$ tends to $\log(1-|\rc|^2)$ in $L^1(\T)$. 
\smallskip

Consider the sequences $\fa_n^{\pm}$, $\fb_n^{\pm}$, $f_n^{\pm} = \fa_n^\pm/\fb_n^\pm$, $a_n$ and $b_n$ corresponding to $q_n$ in a sense described at the beginning of Section \ref{s52}. Furthermore, let $ A$ be an outer function in $\D$ with $ A(0) > 0$ such that  $| A|^{-2} = 1-|\rc|^2$ and $ B = \rc  A$. From the equation \eqref{eq2new}  and the definitions of $a_n$, $b_n$, $A$, $B$, we see that $a_n\to  A$, $b_n\to B$ in Lebesgue measure on $\T$. Also we have $a_n\to  A$ locally uniformly in $\D$.
The functions $\fa_{n}^{\pm}/a_n$ are in the unit ball of $H^2(\D)$ by Proposition \ref{pr5},
hence one can choose a subsequence $n_j$, some functions $\tilde{\fa}^{\pm}$ and Schur functions $\tilde f^{\pm}$ such that
\begin{itemize}
\item $\fa_{n_j}^{\pm} \to \tilde\fa^{\pm}$ locally uniformly in $\D$;
\item $\fa_{n_j}^{\pm}/a_{n_j} \to \tilde\fa^{\pm}/A$ locally uniformly in $\D$ and weakly in $H^2(\D)$;
\item $1/\fa_{n_j}^{\pm} \to 1/\tilde\fa^{\pm}$ locally uniformly in $\D$ and weakly in $H^2(\D)$;
\item $f_{n_j}^{\pm} \to \tilde f^{\pm}$ locally uniformly in $\D$.
\end{itemize}
With this choice of $\tilde\fa^{\pm}$, both functions $\tilde\fa^{\pm}/A$, $A/\tilde\fa^{\pm}$ belong to the Smirnov class in $\D$, hence $\tilde\fa^{\pm}$ are outer functions. Put $\tilde\fb^{\pm} = \tilde f^{\pm} \tilde\fa^{\pm}$. Let $q$ be defined in terms of recurrence coefficients of $\tilde f^{\pm}$ by 
$$q(k) 
= 
\begin{cases}
(\tilde f^+)_{k}(0), & k \ge 0,\\ 
-\ov{(\tilde f^-)_{-k}(0)}, & k < 0.
\end{cases}
$$ 
Note that $(\tilde f^-)_{0}(0) = (\tilde f^-)(0) = 0$ because $(f_n^-)(0) = 0$ for every $n$. We claim that $q_{n_{j}} \to q$ in $\ell^{2}(\Z, \D)$. To prove this, introduce $\fa^{\pm}$, $\fb^{\pm}$, $f^{\pm} = \fa^\pm/\fb^\pm$, $a$, $b$ as the objects from the beginning of Section \ref{s52} corresponding to $q$. It is clear that $f^\pm = \tilde f^\pm$. Let us show that 
\begin{gather*}
	a = A,\quad b = B,\quad \fa^\pm= \tilde\fa^\pm,\quad\fb^\pm =\tilde\fb^\pm. 
\end{gather*}
We have $f^{\pm} = \fb^{\pm}/\fa^{\pm} = \tilde\fb^{\pm}/\tilde\fa^{\pm}$ by construction, and functions $\fa^{\pm}$, $\tilde\fa^{\pm}$ are outer (we do not know, however, that $1-|f^{\pm}|^2 = |\tilde \fa^{\pm}|^{-2}$). Therefore, there are outer functions $s^{\pm}$ such that $\tilde\fa^{\pm} = s^{\pm}\fa^{\pm}$ and $\tilde\fb^{\pm} = s^{\pm}\fb^{\pm}$. It follows that 
\begin{gather*}
	A = \tilde\fa^{+}\tilde\fa^{-} - \tilde\fb^{+}\tilde\fb^{-} = s^{+}s^{-}(\fa^{+}\fa^{-} - \fb^{+}\fb^{-}) = s^{+}s^{-}a.
\end{gather*}
almost everywhere on $\T$ because this relation holds in $\D$. Now write formula \eqref{RH} for $q_{n_j}$ in the form
\begin{gather*}
\begin{pmatrix}
\fa_{n_j}^+ & \fb_{n_j}^{+} \\
\fb_{n_j}^{-} & 
\fa_{n_j}^{-} 
\end{pmatrix} 
=
\begin{pmatrix}
1 & \ov{b_{n_j}} \\
-b_{n_j} & 1
\end{pmatrix}
\overline{
\begin{pmatrix}
\fa_{n_j}^+ & -\fb_{n_j}^{-} \\
-\fb_{n_j}^{+} & \fa_{n_j}^{-}
\end{pmatrix}^{-1}
}
=
\begin{pmatrix}
1 & \ov{b_{n_j}} \\
-b_{n_j} & 1
\end{pmatrix}
\overline{\frac{1}{a_{n_j}}
\begin{pmatrix}
\fa_{n_j}^- & \fb_{n_j}^{-} \\
\fb_{n_j}^{+} & \fa_{n_j}^{+}
\end{pmatrix}.
}
\end{gather*}
Multiplying both sides by $\frac{1}{a_{n_j}}$, we get 
\begin{gather}\label{eq35}
\frac{1}{a_{n_j}}
\begin{pmatrix}
\fa_{n_j}^+ & \fb_{n_j}^{+} \\
\fb_{n_j}^{-} & 
\fa_{n_j}^{-} 
\end{pmatrix} 
=
\begin{pmatrix}
1/a_{n_j} & \ov{b_{n_j}}/a_{n_j} \\
-b_{n_j}/a_{n_j} & 1/a_{n_j}
\end{pmatrix}
\overline{\frac{1}{a_{n_j}}
\begin{pmatrix}
\fa_{n_j}^- & \fb_{n_j}^{-} \\
\fb_{n_j}^{+} & \fa_{n_j}^{+}
\end{pmatrix}.
}
\end{gather}
By construction, we have $\fa_{n_j}^{\pm}/a_{n_j} \to \tilde\fa^{\pm}/A$,
$\fb_{n_j}^{\pm}/a_{n_j} \to \tilde\fb^{\pm}/A$
weakly in $H^2$. We also have  $\ov{b_{n_j}}/a_{n_j} \to \ov{B}/A$, $b_{n_j}/a_{n_j} \to B/A$, $1/a_{n_j} \to 1/A$ strongly in $L^2(\T)$ by the dominated convergence theorem, because $b_{n_j}/a_{n_j}$, $\ov{b_{n_j}}/a_{n_j}$, $1/a_{n_j}$ are uniformly bounded and converge in Lebesgue measure on $\T$ to $B/A$. It follows that both sides of \eqref{eq35} converge weakly in $L^2(\T)$. Taking the limit in \eqref{eq35}, we obtain 
\begin{gather*}
\frac{1}{A}
\begin{pmatrix}
\tilde\fa^+ & \tilde\fb^{+} \\
\tilde\fb^{-} & 
\tilde\fa^{-} 
\end{pmatrix} 
=
\begin{pmatrix}
1/A & \ov{B}/A \\
-B/A & 1/A
\end{pmatrix}
\overline{\frac{1}{A}
\begin{pmatrix}
\tilde\fa^- & \tilde\fb^{-} \\
\tilde\fb^{+} & \tilde\fa^{+}
\end{pmatrix},
}
\end{gather*}
or, in equivalent form,
\begin{gather*}
\begin{pmatrix}
s^{+}\fa^+ & s^{+}\fb^{+} \\
s^{-}\fb^{-} & 
s^{-}\fa^{-} 
\end{pmatrix} 
=
\begin{pmatrix}
1 & \ov{B} \\
-B & 1
\end{pmatrix}
\overline{\frac{1}{s^{+}s^{-}a}
\begin{pmatrix}
s^{-}\fa^- & s^{-}\fb^{-} \\
s^{+}\fb^{+} & s^{+}\fa^{+}
\end{pmatrix}.
}
\end{gather*}
Equation \eqref{RH} written for $q$, $a$, $b$, $\fa^{\pm}$, $\fb^{\pm}$ says
\begin{equation}\notag
\begin{pmatrix}
\fa^+ & \fb^{+} \\
\fb^{-} & 
\fa^{-} 
\end{pmatrix} 
=
\begin{pmatrix}
1 & \ov{b} \\
-b & 1
\end{pmatrix}
\overline{\frac{1}{a}
\begin{pmatrix}
\fa^- & \fb^{-} \\
\fb^{+} & \fa^{+}
\end{pmatrix}.
}
\end{equation}
It follows that 
\begin{align*}
\begin{pmatrix}
s^+ & 0 \\
0 & s^{-} 
\end{pmatrix}
\begin{pmatrix}
1 & \ov{b} \\
-b & 1
\end{pmatrix}
\overline{\frac{1}{a}
\begin{pmatrix}
\fa^- & \fb^{-} \\
\fb^{+} & \fa^{+}
\end{pmatrix}
}
&=\begin{pmatrix}
s^{+}\fa^+ & s^{+}\fb^{+} \\
s^{-}\fb^{-} & 
s^{-}\fa^{-} 
\end{pmatrix} 
=
\begin{pmatrix}
1 & \ov{B} \\
-B & 1
\end{pmatrix}
\overline{\frac{1}{s^{+}s^{-}a}
\begin{pmatrix}
s^{-}\fa^- & s^{-}\fb^{-} \\
s^{+}\fb^{+} & s^{+}\fa^{+}
\end{pmatrix}
}\\
&=
\begin{pmatrix}
1 & \ov{B} \\
-B & 1
\end{pmatrix}
\ov{
\frac{1}{s^{+}s^{-}}
\begin{pmatrix}
s^- & 0 \\
0 & s^{+} 
\end{pmatrix}
}
\overline{\frac{1}{a}
\begin{pmatrix}
\fa^- & \fb^{-} \\
\fb^{+} & \fa^{+}
\end{pmatrix}.
}
\end{align*}
From here we get
\begin{align*}
\begin{pmatrix}
1 & \ov{b} \\
-b & 1
\end{pmatrix}
&=
\begin{pmatrix}
1/s^{+} & 0 \\
0 & 1/s^{-}
\end{pmatrix}
\begin{pmatrix}
1 & \ov{B} \\
-B & 1
\end{pmatrix}
\ov{
\frac{1}{s^{+}s^{-}}
\begin{pmatrix}
s^- & 0 \\
0 & s^{+} 
\end{pmatrix}
}
\\
&=
\begin{pmatrix}
1/s^{+} & 0 \\
0 & 1/s^{-}
\end{pmatrix}
\begin{pmatrix}
1 & \ov{B} \\
-B & 1
\end{pmatrix}
\ov{
\begin{pmatrix}
1/s^+ & 0 \\
0 & 1/s^{-} 
\end{pmatrix}
}\\
&=
\begin{pmatrix}
1/|s^{+}|^2 & \ov{B}/s^{+}\ov{s^{-}} \\
-B/\ov{s^{+}}s^{-} & 1/|s^{-}|^2
\end{pmatrix}.
\end{align*}
It follows that $|s^{\pm}|^2 = 1$. Recall that $s^\pm$ are outer and $s^\pm(0) > 0$, therefore $s^{\pm} = 1$, $\tilde\fa^{\pm} = \fa^{\pm}$, $\tilde\fb^{\pm} = \fb^{\pm}$, $a = A$, $b = B$, and $\rc_{q} = a/b = A/B = \rc$. It remains to show that $q_{n_j} \to q$ in $L^2(\T)$. Since $f^{\pm}$ are locally uniform limits of $f_{n_j}^{\pm}$ in $\D$, we have $\lim_{j \to +\infty} q_{n_{j}}(k) = q(k)$ for each $k \in \Z$ from Schur's algorithm \eqref{eq: NFLT schur algorithm} for $f^{\pm}$. Moreover, \eqref{eq7} and $a = A$ imply
\begin{align*}
\log\prod_{k \in \Z}(1-|q(k)|^2) 
= \log|a(0)|^{-2} = \log|A(0)|^{-2} = \lim_{j \to +\infty}\log|a_{n_j}(0)|^{-2}= \lim_{j \to \infty}\log\prod(1-|q_{n_{j}}(k)|^2).
\end{align*}
The last relation together with elementwise convergence $\lim_{j \to +\infty} q_{n_{j}}(k) = q(k)$ gives $q_{n_{j}} \to q$ in the norm of $\ell^2(\Z, \D)$. \qed 

\begin{Prop}\label{pr10}
The set $\cG = \cup_{\delta \in [0,1)} \cG[\delta]$ is dense in $\ell^2(\Z, \D)$. In fact, $\ell^1(\Z, \D) \subset \cG$. If $q \in \cG$ and $\supp q \subset \Z_+$, then $\|\fc_q\|_{L^\infty(\T)} < 1$ for the function $\fc_{q} = f^{+}$ (see Lemma \ref{l34}). 
\end{Prop}
\beginpf By Baxter's theorem (see Chapter 5 in \cite{Simonbook1}), every measure $\mu$ with recurrence coefficients in $\ell^1(\Z_+, \D)$ has its Szeg\H{o} function, $D_{\mu}$, in the Wiener algebra $W(\T)$. It follows that $\fa^{\pm}$, $\fb^{\pm}$ are continuous and uniformly bounded on $\T$ if $q = \{q(n)\}_{n \in \Z}$ is in $\ell^1(\Z, \D)$, hence the function $a = \fa^+\fa^{-} - \fb^+\fb^-$ is uniformly bounded on $\T$ as well. Formula \eqref{eq2new} then implies that $\rc_q \in B[\delta]$, $q \in \cG[\delta]$, for some $\delta \in [0,1)$. The rest of the proposition is straightforward. \qed

\subsection{Proof of Theorem \ref{t5}.} Recall that the scattering map (or the nonlinear Fourier transform) is defined by 
$$
\F_{sc}: q \mapsto \rc_q,
$$
on the set of sequences $\ell^2(\Z, \D)$, see Proposition \ref{pr1}.  
Assertions $(1)$, $(2)$ of the theorem are Propositions \ref{pr7}, \ref{pr8}, respectively. Assertion  $(3)$ for compactly supported $q: \Z \to \D$ is Proposition~\ref{pr9}. Since $\F_{sc}$ is continuous, assertion $(3)$ then holds for all $q \in \ell^2(\Z,\D)$. To prove assertion $(4)$, consider potentials $q \in \ell^2(\Z, \D)$ supported on $\Z \cap (-\infty, 0)$ and observe that Proposition \ref{pr6} implies $X_+ \subset \F_{sc}(\ell^2(\Z, \D))$. Then, since the set $\F_{sc}(\ell^2(\Z, \D))$ is invariant under multiplication by $z^n$, $n \in \Z$, by assertion $(3)$, the set $\F_{sc}(\ell^2(\Z, \D))$ contains trigonometric polynomials $p$ such that $\|p\|_{L^\infty(\T)}<1$ of arbitrary degree. We claim that the set of such polynomials is dense in $X$. Indeed, one can approximate an arbitrary element of $ X$ by a sequence of continuous functions in the open unit ball of $L^\infty(\T)$ using Lusin's theorem, and then uniformly approximate these continuous functions by Fejer means of their Fourier series. Since 
$\F_{sc}$ is a closed map, the fact that $\F_{sc}(\ell^2(\Z,\D))$ contains a dense subset of $X$ implies $\F_{sc}(\ell^2(\Z,\D)) = X$, and $(4)$ follows. Assertion $(5)$ is Proposition~\ref{pr11}. To prove $(6)$, note that $\F_{sc}(\cG[\delta]) \subset B[\delta]$ by definition and $\F_{sc}(\cG[\delta]) \supset B[\delta]$ because $\F_{sc}: \ell^2(\Z, \D) \to X$ is surjective. Thus, $\F_{sc}: \cG[\delta] \to B[\delta]$ is a continuous surjection. By Proposition \ref{pr4}, this map is injective. Then $\F_{sc}: \cG[\delta] \to B[\delta]$ is a  closed continuous bijection between two topological spaces hence it is a homeomorphism, which is $(6)$. Assertion $(7)$ is not proved in our paper, the reader can find its proof at the end of Chapter 2 in \cite{TTT}. \qed

\section{Appendix}\label{appendix}
Denote by $\ell^0(\Z, \D)$ the set of all sequences $q=\{q_n\}_{n\in \Z}$ such that $|q_n| < 1$ for all $n\in \Z$. In this section we show that for every $q_0\in \ell^0(\Z, \D)$, Ablowitz-Ladik equation \eqref{al} has the unique global solution.
\begin{Lem}[Boundedness, \cite{Gol06}, page 4]\label{lemma: Golinskii}
If $q$ solves \eqref{al} on $[0, t_0]$ for the initial data $q_0\in \ell^0(\Z, \D)$, then $q(t,\cdot)\in \ell^0(\Z, \D)$ for all $t \in [0, t_0]$.
\end{Lem}
\beginpf 
Put $\rho_n(t)^2 = 1 - |q(t,n)|^2$, and assume that for some $n\ge 0$ there exists $t_1\in [0, t_0]$ such that $\rho_n(t_1) = 0$ and $\rho_n(t) > 0$ for all $t\in (0, t_1)$. Then for all $t < t_1$ we have 
\begin{gather*}
	2\rho_n\rho_n' = (\rho_n^2)' = -2\Re\left(q_nq_n'\right) = -2\Re\left(q_n\cdot(i\rho_n^2(q_{n - 1} + q_{n + 1})\right) = 2\rho_n^2\Im\left(q_nq_{n - 1} + q_nq_{n + 1}\right),
	\\	
	\rho_n' = \rho_n\Im\left(q_nq_{n - 1} + q_nq_{n + 1}\right),
	\\
	\rho_n(t) = \rho_n(0)\exp\left[\int_0^t \Im\left(q_nq_{n - 1} + q_nq_{n + 1}\right)\, ds\right].
\end{gather*}
If we now send $t$ to $t_1$, the left hand side will tend to $0$, while the right hand side will not, a contradiction.
\qed

\begin{Lem}[Uniqueness, \cite{Ang2020}\label{l62}, page 20]\label{Lemma: uniqueness}
If $q$, $\tilde q$ solve \eqref{al} on $[0, t]$ for some initial data $q_0\in \ell^0(\Z, \D)$, then $q = \tilde q$.
\end{Lem}
\beginpf
Let $q(t, n)$ and $\tilde{q}(t,n)$ be two solutions for the same initial data $q_0$. We have
\begin{align*}
	-i(q_n' - \tilde{q}_n') &= (1 - |q_n|^2)(q_{n-1} + q_{n + 1}) - (1 - |\tilde{q}_n|^2)(\tilde{q}_{n-1} + \tilde{q}_{n + 1}) 
	\\
	&= (q_{n-1} - \tilde{q}_{n-1}) + (q_{n+1} - \tilde{q}_{n+1}) - (|q_n|^2q_{n-1} - |\tilde{q}_n|^2\tilde{q}_{n-1}) - (|q_n|^2q_{n+1} - |\tilde{q}_n|^2\tilde{q}_{n+1}).
\end{align*}
By Lemma \ref{lemma: Golinskii}, both $|q_n|$ and $|\tilde{q}_n|$ do not exceed $1$ hence
\begin{gather}\label{eq45}
|q_n' - \tilde{q}_n'|\le 2|q_{n - 1} - \tilde{q}_{n - 1}| + 2|q_{n + 1} - \tilde{q}_{n + 1}|  + 4|q_n - \tilde{q}_n|.
\end{gather}
Therefore
\begin{gather*}
(|q_n(t) - \tilde{q}_n(t)|^2)' = 2\Re\left((q_n - \tilde{q}_n)\ov{(q_n' - \tilde{q}_n')}\right),
\\
(|q_n(t) - \tilde{q}_n(t)|^2)'\le 12|q_n - \tilde{q}_n|^2 + 2|q_{n - 1} - \tilde{q}_{n - 1}|^2 + 2|q_{n + 1} - \tilde{q}_{n + 1}|^2
\end{gather*}
Define
\begin{gather*}
	M(t) = \sum_{n\in \Z}\frac{|q_n(t) - \tilde{q}_n(t)|^2}{1 + n^2}.
\end{gather*}
We have $M(0) = 0$ and 
\begin{gather*}
	M'(t) = \sum_{n\in \Z}\frac{(|q_n(t) - \tilde{q}_n(t)|^2)'}{1 + n^2}\le 20 M(t).
\end{gather*}
Then Gr\"onwall inequality gives $M(t) = 0$ for all $t\ge 0$ hence $q$ and $\tilde{q}$ coincide. 
\qed

\begin{Prop}[Existence, \cite{Tbook}, Section 1.1]
For every $q_0\in \ell^0(\Z, \D)$ there exists the unique classical global solution $q$ of \eqref{al}.
\end{Prop}
\beginpf Uniqueness follows from Lemma \ref{l62}.
Rewrite \eqref{al} in the integral form:
\begin{gather}\label{integral al}
	q(t, n) = q_0(n) + \int_0^t i\bigl(1 - |q(s,n)|^2\bigr)\bigl(q(s,n-1) + q(s,n+1)\bigr)\, ds, \quad n \in \Z.
\end{gather}
Equations \eqref{al} and \eqref{integral al} are equivalent. Introduce the space of functions $Y = C([0, t]\times \Z)$ where $t = 1/12$. For $u\in Y$, define the mapping
\begin{gather*}
	F(u)(t, n) = i(1 - |u(t, n)|^2)\bigl(u(t,n-1) + u(t,n+1)\bigr), \quad n\in \Z.
\end{gather*}
In this notation \eqref{integral al} becomes $q(t, n) = q_0(n) + \int_0^t F(q)(s,n)\,ds$. Further, consider
\begin{gather*}
	\Phi(u)(t, n) = q_0(n) + \int_0^t F(u)(s,n)\,ds, \quad n\in \Z.
\end{gather*}
Then solvability of \eqref{integral al} is equivalent to the existence of a fixed point for $\Phi: Y \mapsto Y$. Let us show that $\Phi$ is a contraction acting on the set $B_Y = \{u\in Y\colon \|u\|_Y \le 2\}$. Notice that
\begin{gather*}
|F(u)(s, n)|\le 6\|u\|_Y, \quad s\le t,\quad n\in \Z,
\\
|\Phi(u)(t, n)|\le |q_0(n)| + \int_0^t|F(u)(s,n)|\, ds\le 1 + 6t\|u\|_Y,\\
\|\Phi(u)\|_Y \le 1 + 6t\|u\|_Y.
\end{gather*}
In particular, if $u\in Y$, then $\Phi(u)\in Y$. Furthermore, from \eqref{eq45} we see that for $u,v\in Y$  we have
\begin{align*}
|\Phi(u)(t,n) - \Phi(v)(t, n)|&\le \int_0^t|F(u)(s, n) - F(v)(s, n)|\, ds \le 6t\|u - v\|_Y.
\end{align*}
We have $6t < 1$, hence $\Phi$ is a contraction and \eqref{al} has a solution on $[0, t]$.  By Lemma \ref{lemma: Golinskii}, $q(t, \cdot)$ also satisfies $q(t, n)< 1$ for all $n\in \Z$, hence the fixed point algorithm can be applied to find the solution on the segment $[t,2t]$. Iterating this procedure, we obtain the existence of a solution on $[0, \infty)$. The similar argument works for negative $t$, hence the proof is concluded. \qed 


\medskip

The following proposition gives a proof of the convergence in Theorem  \ref{t2} based on the idea from Lemma \ref{Lemma: uniqueness}.

\begin{Prop}
Take $q_0\in \ell^0(\Z, \D)$ and let $q_{0,N}$, $q$, $q_N$ be as in Theorem \ref{t2}. 
Then, for $N\ge |j|$, $t>0$ and all $r \in (0,1)$, we have
\begin{gather*}
	|q(t, j) - q_N(t,j)|\le \frac{\sqrt{2}re^{10t/r^2}}{\sqrt{1 - r^2}} r^{N-|j|}.
\end{gather*}
If we assume $\ell^2(\Z, \D)$, then 
\begin{gather*}
	|q(t, j) - q_N(t,j)|\le re^{10t/r^2}\sqrt{\sum_{|m| > N}|q_0(m)|^2} \cdot r^{N-|j|}.
\end{gather*}
\end{Prop}
\beginpf Set
$M_N(t) = \sum_{m\in \Z}|q(t, m) - q_N(t, m)|^2r^{2|m|}$. 
At $t = 0$ we have
\begin{gather}\label{eq: ineq at time t = 0}
	M_N(0) = \sum_{|m| > N}|q_0(m)|^2r^{2|m|} \le \sum_{|m| > N}r^{2|m|} = \frac{2r^{2N + 2}}{1 - r^2}.
\end{gather}
The inequalities similar  to \eqref{eq45} give us $M_N'(t)\le 20r^{-2}M_N(t)$, hence
\begin{gather*}
	|q(t, j) - q_N(t, j)|^2r^{2|j|}\le M_N(t)\le \exp(20r^{-2}t)M_N(0) = \frac{2e^{20t/r^2}r^{2N + 2}}{1 - r^2}.
\end{gather*}
The first part of the proposition follows. To establish the second inequality, we change the bound \eqref{eq: ineq at time t = 0}. We have
\begin{gather*}
	M_N(0) = \sum_{|m| > N}|q_0(m)|^2r^{2|m|} \le r^{2(N + 1)}\sum_{|m| > N}|q_0(m)|^2.
\end{gather*}
Therefore
\begin{gather*}
	|q(t, j) - q_N(t, j)|^2r^{2|j|}\le M_N(t)\le \exp(20r^{-2}t)M_N(0) = e^{20t/r^2}r^{2N + 2}\sum_{|m| > N}|q_0(m)|^2,
\end{gather*}
which concludes the proof.
\qed
	
\bibliographystyle{plain} 	\bibliography{bibfile}
\end{document}